\documentclass[a4paper,reqno,11pt]{amsart}
\usepackage{amsmath}       
\usepackage[margin=2.5cm]{geometry}
\usepackage[foot]{amsaddr}
\usepackage{amsfonts}      
\usepackage{amssymb}       
\usepackage{amsthm}        
\usepackage{bm}            
\usepackage{enumitem}
\usepackage{graphicx}      
\usepackage{mathtools}     
\usepackage{xcolor}        
\usepackage{dsfont}
\usepackage{tikz}
\usetikzlibrary{positioning,arrows.meta,calc}
\usepackage{caption}

\usepackage[
    colorlinks=true,
    linkcolor=blue,
    citecolor=orange,
    urlcolor=urlblue
]{hyperref} 

\usepackage{aliascnt}
\usepackage[capitalize]{cleveref}

\usepackage{glossaries}
\expandafter\def\csname ver@etex.sty\endcsname{3000/12/31}

\usepackage{autonum}

\newcommand{\R}{\mathbb{R}}
\newcommand{\N}{\mathbb{N}}
\newcommand{\Pp}{\mathcal{P}}
\newcommand{\ent}{\mathrm{Ent}}
\newcommand{\Ric}{\mathrm{Ric}}
\newcommand{\Sec}{\mathrm{Sec}}
\newcommand{\II}{\mathrm{II}}
\newcommand{\Rm}{\mathrm{Rm}}
\newcommand{\Hess}{\mathrm{Hess}}
\newcommand{\Cut}{\mathrm{Cut}}

\newcommand{\spt}{\mathrm{spt}}
\newcommand{\vol}{\mathrm{vol}}

\newcommand{\rcd}{\mathrm{RCD}}

\newcommand{\op}{\mathrm{op}}
\newcommand{\HS}{\mathrm{HS}}
\newcommand{\m}{\mathfrak{m}}
\renewcommand{\S}{\mathbb{S}}
\renewcommand{\d}{\,\mathrm{d}}
\renewcommand{\liminf}{\varliminf}
\renewcommand{\limsup}{\varlimsup}

\usepackage{xpatch} 
\xpatchcmd{\proof} 
{\itshape}
{\bfseries}
{}
{}

\usepackage[square,numbers]{natbib}

\numberwithin{equation}{section}

\newtheorem{theorem}{Theorem}[section]
\newcommand{\newsharedtheorem}[2]{%
  \newaliascnt{#1}{theorem}%
  \newtheorem{#1}[#1]{#2}%
  \aliascntresetthe{#1}}
\newsharedtheorem{lemma}{Lemma}
\newsharedtheorem{corollary}{Corollary}
\newsharedtheorem{proposition}{Proposition}

\theoremstyle{definition}
\newsharedtheorem{definition}{Definition}
\newsharedtheorem{remark}{Remark}

\newsharedtheorem{example}{Example}

\newsharedtheorem{conjecture}{Conjecture}

\makeatletter
\renewcommand\subsubsection{\@startsection{subsubsection}{3}{\z@}%
	{.5\linespacing\@plus.7\linespacing}{-.5em}
	{\normalfont\bfseries}}

\begin{document}
	
\title[Entropic Ricci bounds at conic singularities and boundaries]{Entropic Ricci bounds at conic singularities and boundaries}
\author{Zhenhao Li}

\address{Faculty of Mathematics \\ 
    Bielefeld University\\
	Postfach 10 01 31 \\
	33501 Bielefeld \\
	Germany.}
\email{zhenhao.li@math.uni-bielefeld.de}


\begin{abstract}
	We study synthetic Ricci bounds and related rigidity phenomena in conic and boundary neighborhoods through the behavior of entropy along optimal transports.
	For cones over manifolds, we show that two-sided entropic Ricci bounds at the vertex force the fiber to have constant sectional curvature $1$, while non-trivial quotient singularities, such as Euclidean cones over spherical space forms, remain admissible. An analogous rigidity holds for Riemannian manifolds with conic singularities.
	We also give a local characterization of boundary convexity in terms of entropic lower Ricci bounds, based on a direct optimal transport argument.
\end{abstract}

\maketitle
\tableofcontents

\section{Introduction}

Lower bounds on Ricci curvature admit an intrinsic description via optimal transport. 
On a smooth Riemannian manifold, such bounds are equivalent to suitable convexity properties of entropy along Wasserstein geodesics; see \cite{vonRenesseSturm}.
This formulation is \emph{synthetic} i.e. it does not depend on the existence of an underlying smooth structure.
The pioneering work of Lott--Villani \cite{LottVillani09} and Sturm \cite{sturm2006--1,Sturm06-2} used this principle to define the curvature-dimension condition (CD condition for short) for general metric measure spaces.
Its Riemannian refinement, the $\rcd$ condition introduced by Gigli \cite{GigliRCD}, provides a non-smooth framework, linking the entropic approach via optimal transport with Bakry--\'Emery $\Gamma$-calculus, and various geometric/functional inequalities; see e.g. \cite{AGSRCD2014,AGS15,EKS15}.

The synthetic characterization of upper Ricci bounds is more delicate.
Requiring the entropy to be concave along every optimal transport would already fail in Euclidean space.
Instead, one asks for suitable localized transports along which the entropy is almost concave.
Sturm \cite{Sturm2017remarks} developed synthetic upper Ricci bounds using both this entropic approach and the short-time behavior of the heat flow. 
Their dynamic counterparts lead to the weak and rough notions of Ricci flow respectively, studied in \cite{Erbar-Li-Schultz26}.
On smooth closed manifolds, the entropic formulation for lower/upper Ricci bounds can be summarized heuristically as
\begin{equation}\label{eq:intro-heuristics}
	\begin{aligned}
		\Ric\geq 0
		&\quad\Longleftrightarrow\quad\text{entropy is convex along \emph{every} optimal transport},\\
		\Ric\leq 0
		&\quad\Longleftrightarrow\quad
		\begin{gathered}
			\text{in \emph{every} direction, there \emph{exists} a suitable local transport}\\[-2pt]
			\text{along which entropy is \emph{almost} concave}.
		\end{gathered}
	\end{aligned}
\end{equation}
Here ``almost'' allows an arbitrarily small error relative to the squared transport distance after localization.

 Let $(X,d,\m)$ be a metric measure space, let $W_2$ denote its quadratic Wasserstein distance, and set $\ent_{\m}(\mu)\coloneqq\int\rho\log\rho\d\m$ if $\mu=\rho\m$, and $+\infty$ otherwise.
 Following \cite{Sturm2017remarks,Erbar-Li-Schultz26}, we study the entropic Ricci bounds through the following pointwise quantities.
 \begin{definition}\label{def:eta}
 	For $x,y\in X$ and $\varepsilon>0$, define
 	\begin{equation}\label{eq:intro-eta-epsilon}
 		\eta^\pm_\varepsilon(x,y)\coloneqq\inf_{(\mu^a)}
 		\frac{1}{W_2^2(\mu^0,\mu^1)}\left[ 
 		\partial^\pm_a\ent_{\mathfrak m}(\mu^a)|_{a=1-}
 		-\partial^\mp_a\ent_{\mathfrak m}(\mu^a)|_{a=0+}\right],
 	\end{equation}
 where the infimum is taken over all non-constant $W_2$-geodesics $(\mu^a)_{a\in[0,1]}$ with $\ent_\m(\mu^0)$, $\ent_\m(\mu^1)$ finite, $\spt(\mu^0)\subset B_\varepsilon(x)$, and $\spt(\mu^1)\subset B_\varepsilon(y)$.
 
 We define the lower and upper quantity $\eta_*$ and $\eta^*$ at $p\in X$ as follows
 	\begin{equation}\label{eq:intro-eta}
 		\eta^\pm(x,y)\coloneqq\lim_{\varepsilon\downarrow0}\eta^\pm_\varepsilon(x,y),\qquad \eta_*(p)\coloneqq \liminf_{x,y\to p}\eta^-(x,y),\qquad 
 		\eta^*(p)\coloneqq\limsup_{x,y\to p}\eta^+(x,y).
 	\end{equation}
 \end{definition}
 Here, $\partial^\pm_a\big|_{a=c\pm}$ denotes the upper/lower left/right derivative, and we use the convention $+\infty-\infty=-\infty$.
 
These quantities are consistent with smooth Ricci bounds: at a smooth interior point of a Riemannian manifold $(M,g)$ endowed with $\m\coloneqq e^{-f}\vol_g$, $\eta_*$ and $\eta^*$ recover the smallest and largest eigenvalues of the Bakry-\'Emery Ricci tensor 
 \begin{equation}
 	\Ric_f\coloneqq \Ric_g+\Hess f.
 \end{equation}
 The definition of the $\eta$-quantities nevertheless requires neither a smooth structure nor a pointwise curvature tensor. 
 The question is which features of a singularity they detect.

\subsubsection*{Cones and conic singularities}
Cones provide a natural setting for this question.
Bacher--Sturm \cite{BacherSturm14} established CD bounds for Euclidean and spherical cones over suitable smooth fibers.
For general metric measure spaces, Ketterer \cite{Ketterer-JMPA15} characterized the $\rcd$ condition of cones in terms of the corresponding condition on their fibers.

Adding an upper Ricci bound leads to rigidity: Erbar--Sturm \cite{Erbar-Sturm} proved that a Euclidean cone with $\rcd$ condition admitting a finite upper Ricci bound at the vertex in the heat-flow formulation must be Euclidean space; see also \cite[Appendix~C]{Erbar-Li-Schultz26} for the corresponding rigidity for spherical cones.
Here the upper bound is defined through the short-time asymptotics of the $W_2$-distance between heat kernels.
As one of the main contributions of this work, we show that the entropic formulation has a different rigidity threshold.

Let $(F,h,e^{-f}\vol_h)$ be a compact smooth weighted Riemannian manifold of dimension $n\geq2$, possibly with smooth boundary.
For $K\in\mathbb R$, $N\geq n$, denote by $\Sigma_K^N(F)$ the $(K,N)$-cone over the metric measure space $(F,h,e^{-f}\vol_h)$; see \cref{sec:cones} for the precise definition.
\begin{theorem}[Cone rigidity]\label{thm:intro-exact-cone}
	Let $\Sigma_K^N(F)$ be the $(K,N)$-cone of $(F,h,e^{-f}\vol_h)$ and let $o$ be a vertex of $\Sigma_K^N(F)$. The following are equivalent:
	\begin{enumerate}
		\item $\eta_*\equiv\eta^*\equiv nK$ on $\Sigma_K^N(F)$;
		\item $\eta_*(o)>-\infty$ and $\eta^*(o)<+\infty$;
		\item $N=n$, $f$ is constant, $\Sec_h\equiv1$, and $\II\geq0$ on $\partial F$.
	\end{enumerate}
\end{theorem}
Here $\II$ denotes the second fundamental form of $\partial F$, defined by $\II(X,Y)\coloneqq-h(\nabla_X\nu,Y)$, where $\nu$ is the inward unit normal vector; see \cref{sec:boundary-rigidity}.
We postpone the discussion of boundary rigidity until later in the introduction.

Finiteness of the two local quantities at a single vertex thus forces the fiber to be unweighted, with constant sectional curvature $1$ and convex boundary.
The vertex may nevertheless remain singular: for example, the Euclidean cone over the round $\mathbb{RP}^n$, equipped with Riemannian volume, satisfies $\eta_*\equiv\eta^*\equiv0$.
The distinction may be expressed as follows:
\begin{center}
	\emph{The entropic upper bound detects curvature blow-up near the vertex}\par
	\emph{without necessarily detecting the vertex as a singular point.}
\end{center}
This is compatible with the fact that, under the curvature-dimension assumptions of \cite{BacherSturm14}, almost every path of an optimal transport between absolutely continuous measures avoids the vertex. The obstruction instead arises from the geometry of transports in the regular part.

The mechanism can be seen through the Jacobian calculus.
Let $(\mu^a)$ be a (suitably regular) $W_2$-geodesic and denote by $\gamma_x$ the transport geodesic starting at $x$.
Then the entropy satisfies
\begin{equation}\label{eq:intro-jacobian}
	\frac{\d^2}{\d a^2}\ent_{\m}(\mu^a)
	=\int_M|S_x(a)|_{\HS}^2+\Ric_f(\dot\gamma^a_x,\dot\gamma^a_x)\d\mu^0(x).
\end{equation}
Here $S_x$ solves the matrix Riccati equation $S_x'+S_x^2+R_x=0$ along $\gamma_x$, where $R_x(a)\coloneqq\Rm(\cdot,\dot\gamma^a_x)\dot\gamma^a_x$ and $\Rm$ denotes the Riemann curvature tensor. 
If $R_x$ does not vanish along $\gamma_x$, the Riccati equation yields a positive lower bound for $\int_0^1|S_x|_{\HS}^2\d a$, uniformly over nearby geodesics and independently of the initial data for $S_x$. 
On the regular part of a cone with non-negative Ricci, this produces a strictly positive entropic defect for all (regular) transports between sufficiently small neighborhoods of given endpoints. 
Contracting this configuration towards the vertex by a dilation of factor $\lambda\searrow0$ amplifies the defect by $\lambda^{-2}$, forcing $\eta^*(o)=+\infty$; see \cref{sec:coercivity-Riem,thm:rigidity_E_cone}.

The usual derivation of displacement convexity from a lower Ricci bound discards the non-negative term $|S_x|_{\HS}^2$.
Our argument exploits the dependence of this term on the full Riemannian curvature through the Riccati equation.
A related matrix perspective was introduced by Shenfeld \cite{Shenfeld24} and developed on Riemannian manifolds by Aishwarya--Rotem--Shenfeld \cite{ARS26}, who characterized non-negative sectional curvature through matrix displacement convexity of a tensorial refinement of entropy.
Our cone rigidity result illustrates how the scalar entropy can already impose further geometric constraints.

Exact homogeneity under dilations is unavailable when $K\neq0$.
We therefore establish a rigidity result for asymptotically conic neighborhoods, which covers exact cones for general $K\neq 0$. 
Near a conic singularity $o$, a punctured neighborhood can be identified with $(0,l)\times F$, where $l>0$ and $F$ is a smooth manifold without boundary.
In this coordinate, the metric takes the form $g=\d r^2+r^2h_r$, where $(h_r)_{r\in(0,l)}$ is a smooth family of Riemannian metrics on $F$ converging to a Riemannian metric $h_0$ as $r\to0$; see \cref{sec:conicnghd} for the precise definition.
The corresponding limiting cone has metric $g_\Sigma=\d r^2+r^2h_0$ over the fiber $(F,h_0)$.
Under the additional asymptotic condition on the metric, namely,
\begin{equation}\label{eq:intro-conic-decay}
	\big| (\nabla^{g_\Sigma})^j(g-g_\Sigma)\big|_{g_\Sigma}
	=o(r^{\alpha-j}),\qquad j=0,1,2;
\end{equation}
for some $\alpha>0$, we can transfer the entropic defect from the limiting cone to the original conic neighborhood.

	\begin{theorem}[Rigidity of conic singularities]\label{thm:intro-conic}
		Let $(M,d,g,o)$ be an $(n+1)$-dimensional Riemannian manifold with a conic singularity of order $2$ in the sense of \cref{def:conic-singularity}, where $n\geq2$.
       If
		\begin{equation}\label{eq:intro-tip-bounds}
			\eta_*(o)>-\infty,\qquad \eta^*(o)<+\infty,
		\end{equation}
		then $\Sec_{h_0}\equiv1$, where $(F,h_0)$ is the fiber of the limiting cone. 
		Under these assumptions, the following are equivalent:
		\begin{enumerate}
			\renewcommand{\labelenumi}{(\alph{enumi})}
			\item $o$ is a topological manifold point of $M$;
			\item $(F,h_0)$ is isometric to the unit round sphere $\mathbb S^{n}$;
			\item the volume density satisfies
			\begin{equation}\label{eq:intro-density}
				\Theta(o)\coloneqq\lim_{r\searrow0}
				\frac{\vol_g(B_r(o))}{\omega_{n+1} r^{n+1}}=1.
			\end{equation}
		\end{enumerate}
	\end{theorem}
In dimension two, a different phenomenon occurs: for a Euclidean cone of angle $0<\theta<2\pi$, $\eta_*=\eta^*=0$ everywhere; see \cite[Proposition 6.2]{Erbar-Li-Schultz26}.

\subsubsection*{Boundary rigidity/non-rigidity}

A non-convex boundary can obstruct lower curvature bounds in the sense of Alexandrov, even for smooth Riemannian manifolds of bounded geometry, since geodesics may branch at the boundary.
A similar obstruction occurs for synthetic lower Ricci bounds: via Gigli's measure-valued Ricci tensor, Han \cite{Han-17,Han-20} showed for smooth manifolds with boundary that RCD conditions are equivalent to the corresponding lower bound on Ricci curvature together with the convexity of the boundary.

The situation is different for upper curvature bounds.
Every smooth boundary point admits a neighborhood with curvature bounded above in the sense of Alexandrov; see \cite{ABB1993}.
Using the semiconvexity for squared distances between geodesics in such neighborhoods, we obtain a quantitative upper bound for $\eta^*$ at smooth boundary points; see \cref{sec:CAT-upper}.
Together with our local characterization of boundary convexity via $\eta_*$, which also enters the proof of \cref{thm:intro-exact-cone}, this yields the following statement.

\begin{theorem}\label{thm:intro-boundary}
	Let $(M,g,e^{-f}\vol_g)$ be a smooth weighted Riemannian manifold with smooth boundary, and let $p\in\partial M$.
	Then $\eta^*(p)<+\infty$, and
	\begin{equation}\label{eq:intro-boundary}
		\eta_*(p)>-\infty
		\quad\Longleftrightarrow\quad
		\II\geq0\ \text{on a neighborhood of }p\text{ in }\partial M.
	\end{equation}
\end{theorem}

The lower-bound characterization is proved by a direct transport argument: a negative direction of $\II$ produces optimal transports that concentrate mass on the boundary, giving rise to arbitrarily negative entropy defects; see \cref{sec:boundary-rigidity}.

\subsubsection*{Comparison with the synthetic Ricci bounds via heat flow}
For comparison, we recall another synthetic Ricci upper bound via heat flow mentioned earlier, introduced by Sturm \cite{Sturm2017remarks}.
Denote by $\widehat P_t$ the dual heat flow acting on probability measures and, for $x\neq y$, define
\begin{equation}\label{eq:intro-theta}
	\vartheta^\pm(x,y)\coloneqq -\partial^\mp_t \log \left(W_2(\widehat P_t\delta_x,\widehat P_t\delta_y)\right)|_{t=0},
	\quad
	\vartheta^*(p)\coloneqq\limsup_{x,y\to p}\vartheta^+(x,y),
	\quad
	\vartheta_*(p)\coloneqq\liminf_{x,y\to p}\vartheta^-(x,y).
\end{equation}
On smooth closed manifolds, $\vartheta^*$ and $\vartheta_*$ recover the largest and smallest eigenvalues of the Ricci tensor, respectively; see \cite[Theorem 3.1]{Sturm2017remarks}. 
At cone vertices, finiteness of $\vartheta^*(o)$ additionally rules out the non-trivial Euclidean quotients allowed by \cref{thm:intro-exact-cone}; see \cite{Erbar-Sturm}. 
The Euclidean halfspace gives another distinction regarding boundary rigidity: $\eta_*=\eta^*=0$ at its boundary, whereas one expects  $\vartheta^*=+\infty$ there, as suggested by regarding the halfspace as a cone over a hemisphere and applying \cite[Theorem~1.1]{Erbar-Sturm}; a rigorous verification of $\vartheta^*=+\infty$ at boundary points, however, remains open.

The distinction extends to synthetic Ricci flows. 
Let $(F,h)$ be a closed $n$-dimensional manifold, $n\geq2$, with $\Ric_h=(n-1)h$. On its spherical cone $\Sigma$, consider
\begin{equation}\label{eq:intro-flow}
	d_t=(1-2nt)^{1/2}d_{\Sigma},\qquad
	\mathfrak m_t=(1-2nt)^{(n+1)/2}\vol_{\Sigma},
	\qquad 0\leq t<\frac1{2n}.
\end{equation}
As an application of \cref{thm:intro-exact-cone}, we show that this time-dependent metric measure space is a weak Ricci flow exactly when $\Sec_h\equiv1$, and a rough Ricci flow only when $F$ is the unit round sphere; both in the sense of \cite{Erbar-Li-Schultz26}. 
Non-trivial spherical space forms therefore separate the two notions in higher dimensions, extending the two-dimensional comparison in \cite[Proposition~6.2]{Erbar-Li-Schultz26}.
Another application is the instability of the entropic synthetic Ricci upper bound under the Gromov--Hausdorff convergence; cf. \cite[Section~5]{Sturm2017remarks} for a stable notion of synthetic upper Ricci bound.
We refer the reader to \cref{sec:application} for details.

These examples suggest a stronger rigidity principle regarding the $\vartheta$-quantities. 
We conjecture, for a manifold with a conic singularity $o$ as in \cref{thm:intro-conic}, under an appropriate regularity assumption, that 
\begin{equation}
	\vartheta^*(o)<+\infty
	\quad\Longrightarrow\quad	
	(F,h_0)\cong\mathbb S^{n}.
\end{equation}
\subsubsection*{Notation and Convention}
Unless otherwise specified, all manifolds are smooth and connected, possibly with smooth boundary and not necessarily complete. They are equipped with a smooth Riemannian metric and the induced intrinsic distance.

The letters $g,h$ are reserved for Riemannian metrics, and ${\rm L}_g$, $\vol_g$, $\Ric_g$, $\Hess_g$ denote the length of curves, the Riemannian volume measure, Ricci tensor and the Hessian with respect to $g$. 
The letter $f$ is reserved for smooth weight functions. Thus $\Ric_f$ always denotes the weighted Ricci tensor and should not be confused with $\Ric_g$.

Given a reference measure $\m$, the relative entropy of $\mu$ is $\ent_{\m}(\mu)\coloneqq\int\rho\log\rho\d\m$ if $\mu=\rho\m$, and $+\infty$ otherwise.
Unless a weight is specified, we consider the Riemannian volume measure and in that case we simply write $\ent$ without any subscript.

Finally, $\mathcal{L}^n$ denotes the $n$-dimensional Lebesgue measure, and $\omega_n$ denotes the volume of the $n$-dimensional unit ball. For a linear map between finite-dimensional inner product spaces, $|{\cdot}|_{\HS}$ and $|{\cdot}|_{\mathrm{op}}$ denote its Hilbert--Schmidt norm and its operator norm.

%
\subsubsection*{Use of AI tools.}
All theorems and proofs were originally developed by the author.
ChatGPT and later Claude were then used to check the mathematics (the proof of \cref{cor:finite-at-boundary} was revised accordingly), to polish the language, and to generate the figures.
Besides, the Stenzel metric example in \cref{sec:application} was found by ChatGPT, which also assisted substantially in drafting that part.
\section{Preliminaries}
\subsection{Synthetic Ricci bounds}\label{sec:eta_pre}
	The curvature-dimension condition $\mathrm{CD}(K,N)$ encodes a lower
	Ricci curvature bound $K$ and an upper dimension bound $N$ through
	convexity properties of entropy along optimal transports; see
	\cite{LottVillani09,sturm2006--1,Sturm06-2}.
	Its Riemannian refinement, the $\rcd(K,N)$ condition, additionally
	requires the Cheeger energy to be quadratic.
	We refer to \cite{GigliRCD,AGSRCD2014,AGS15,EKS15} for the precise
	definitions and equivalent formulations.
	
	We use the entropic quantities $\eta^\pm_\varepsilon$, $\eta^\pm$, $\eta_*$ and $\eta^*$ introduced in \cref{def:eta}.
	Their consistency with classical Ricci bounds follows from \cite{Sturm2017remarks} as well as from the static specialization of \cite[Corollary~4.12]{Erbar-Li-Schultz26}.
	More precisely, let $(M,g,e^{-f}\vol_g)$ be a smooth closed weighted Riemannian manifold, and write $\Ric_f\coloneqq\Ric_g+\operatorname{Hess}_g f$.
	Then, for every $p\in M$,
	\begin{equation}\label{eq:eta-smooth-consistency}
		\eta_*(p)=\min_{\substack{v\in T_pM\\ |v|_g=1}}\Ric_f(v,v),
		\qquad
		\eta^*(p)=\max_{\substack{v\in T_pM\\ |v|_g=1}}\Ric_f(v,v).
	\end{equation}
	By localization, the identities in
	\eqref{eq:eta-smooth-consistency} also hold at smooth interior points
	of manifolds possibly incomplete or with boundary.
	
	On an $\rcd(K,\infty)$ space, the entropy is $K$-convex along every $W_2$-geodesic with finite-entropy endpoints.
    Then 
	\begin{equation}\label{eq:eta-RCD-lower}
		\eta^+_\varepsilon(x,y)=\eta^-_\varepsilon(x,y)\geq K
		\end{equation}
	for all $\varepsilon>0$ and $x,y$; see \cite[Section~2.3 and Lemma~3.5]{Erbar-Li-Schultz26} or \cite[Section~2.2.3]{Li26}.
	More generally, \eqref{eq:eta-RCD-lower} also holds on essentially non-branching $\mathrm{CD}(K,N)$ spaces with $N<\infty$.
	Equality of the two sign choices follows whenever the entropy is upper regular in the sense of \cite{Sturm2018Super}; see also \cite[Definition~2.3]{Erbar-Li-Schultz26}.
	These hypotheses need not hold for general manifolds with boundary.
	
	Nevertheless, the choice of signs does not affect the estimates in this paper.
	Lower estimates on $\eta$ are established for the smaller one $\eta^-$.
	Upper estimates are obtained using explicit transports for which the entropy admits derivatives at both endpoints.
	We therefore do \textbf{not} distinguish between $\eta^+$ and $\eta^-$ in what follows.
\subsection{Jacobian calculus along optimal transports}\label{sec:Jacobian-calculus}
The term ``Jacobian calculus" refers here to the use of Jacobi fields to estimate the Jacobian determinants of optimal transport maps on Riemannian manifolds. 
This is a standard tool in the study of displacement convexity; see e.g. \cite{Cordero-Erausquin2001,Sturm05,vonRenesseSturm,LottVillani09}. 
In this subsection, we recall the basic ideas and fix the notations used in the subsequent sections.

Let $(M,g)$ be a smooth complete Riemannian manifold without boundary and let $\m\coloneqq e^{-f}\vol_g$ be the reference measure where $f$ is a smooth function on $M$.
Let $\mu^0,\mu^1\in\Pp_2(M)$ with compact supports and assume that $\mu^0,\mu^1\ll\vol_g$.
By \cite{Mccann2001}, the unique optimal coupling is induced by a map
\begin{equation}
	T(x)=\exp_x(-\nabla\varphi(x)),
\end{equation}
where $\varphi$ is a $d_g^2/2$-concave function, called a
Kantorovich potential. 
The unique $W_2$-geodesic connecting $\mu^0$ and $\mu^1$ is given by
\begin{equation}
	\mu^a=(T^a)_\#\mu^0,\qquad T^a(x)=\exp_x(-a\nabla\varphi(x)),\qquad a\in[0,1].
\end{equation}
For $\mu^0$-a.e. $x$, the Hessian of $\varphi$ exists, and the approximate differential
\begin{equation}\label{eq:Jacobi-matrix}
	 A_x(a)\coloneqq \d T^a(x)\colon T_xM\to T_{T^a(x)}M
\end{equation}
is defined and non-singular for all $a\in[0,1]$; see \cite[Theorem 4.2]{Cordero-Erausquin2001}.

For such $x$, set $\gamma_x(a)\coloneqq T^a(x)$.
For every $v\in T_xM$, the vector field $a\mapsto A_x(a)v$ is a Jacobi field along $\gamma_x$.
Thus $ A_x$ is the unique solution of
\begin{equation}\label{eq:Jacobi-equation}
	\nabla_a\nabla_a A_x+\Rm( A_x,\dot\gamma_x)\dot\gamma_x=0,
	\qquad
	 A_x(0)={\rm Id},\quad\nabla_a A_x(0)=-\Hess\varphi(x).
\end{equation}
Here $\nabla_a$ denotes the covariant derivative along $\gamma_x$, and the Hessian is viewed as an endomorphism via $g$.
Wherever $ A_x(a)$ is non-singular, define
\begin{equation}\label{eq:tensorS}
	 S_x(a)\coloneqq
	\nabla_a A_x(a)\circ A_x(a)^{-1}
	\colon T_{\gamma_x(a)}M\to T_{\gamma_x(a)}M,
\end{equation}
which solves the Riccati type equation
\begin{equation}\label{eq:Riccati-equation}
	\nabla_a S_x+ S_x^2+\Rm(\cdot,\dot\gamma_x)\dot\gamma_x=0.
\end{equation}

\begin{lemma}\label{lemma:Jacobiancalculus}
	Let $\mu^0,\mu^1\in\Pp_2(M)$ have compact supports and finite entropy.
	Then, along the $W_2$-geodesic $(\mu^a)_{a\in[0,1]}$ connecting $\mu^0$ and $\mu^1$,
	\begin{equation}
	\partial^{-}_a \ent_\m(\mu^a)\big|_{a=1-}-\partial_a^{+} \ent_\m(\mu^a)\big|_{a=0+}= \int_M\int_0^1 | S_x(a)|_{\HS}^2+ \Ric_f(\dot\gamma_{x}(a),\dot\gamma_{x}(a)) \d a\d\mu^0(x)\,
\end{equation}
	where $\gamma_x(a)=T^a(x)$ and $ S_x$ is defined by
	\eqref{eq:tensorS}.
\end{lemma}
\begin{proof}
	The proof is standard by the preceding Jacobian calculus and the change-of-variable formula \cite[Corollary~4.7]{Cordero-Erausquin2001} (see also formula (4.16) in \cite{sturm2006--1} for the weighted case).
	We provide a sketch for the reader's convenience.
	Set
		\begin{equation}
			J_x(a)\coloneqq
			f(\gamma_x(a))-f(x)-\log\det A_x(a).
		\end{equation}
		By the change-of-variable formula,
		\begin{equation}
			\ent_\m(\mu^a)=\ent_\m(\mu^0)+\int_M J_x(a)\d\mu^0(x).
		\end{equation}
	Then standard computation gives
	\begin{equation}
		J_x''(a)=\Hess f(\dot\gamma_x(a),\dot\gamma_x(a))-\partial^2_a \log\det A_x(a),\qquad \partial^2_a \log\det A_x(a)={\rm tr}(\nabla_a  S_x(a)),
	\end{equation}
    where $S_x$ solves the Riccati type equation \eqref{eq:Riccati-equation} for $\mu^0$-a.e. $x$.
	Since $ S_x$ is self-adjoint, taking the trace in \eqref{eq:Riccati-equation} (with $\Ric_f=\Hess f+\Ric$) gives
		\begin{equation}
			J_x''(a)={\rm tr}(S^2_x)+\Ric_f(\dot\gamma_x(a),\dot\gamma_x(a))= | S_x(a)|_{\HS}^2+\Ric_f(\dot\gamma_x(a),\dot\gamma_x(a)).
		\end{equation}
	Passing the derivatives through the integral yields the assertion.
\end{proof}
\subsection{Cones over metric measure spaces}\label{sec:cones}
\begin{definition}\label{def:cone}
    Let $(F,d,\m)$ be a metric measure space.
    For $K\in\R$ and $N\geq 1$, the $(K,N)$-cone over $F$ is the metric measure space
    \begin{equation}
        \bigl(\Sigma_K(F),d_{\Sigma_K},\m_{\Sigma_K^N}\bigr),
        \qquad
        \Sigma_K(F)\coloneqq(I_K\times F)/{\sim},
    \end{equation}
    where $I_K\coloneqq[0,\frac{\pi}{\sqrt{K^+}}]$ with $K^+\coloneqq K\vee 0$ and the convention $[0,\frac{\pi}{0}]=[0,\infty)$; and $(r,x)\sim (r,y)$ if $r\in\partial I_K$.
    
    For $(r,x),(r',x')\in\Sigma_K(F)$, set
    \begin{equation}
        d_{\Sigma_K}\bigl((r,x),(r',x')\bigr)
        \coloneqq
        \begin{cases}
            \cos_K^{-1}\bigl(\cos_K(r)\cos_K(r')
            +K\sin_K(r)\sin_K(r')\cos(d(x,x')\wedge\pi)\bigr),
            & K\neq0,\\
            \sqrt{r^2+r'^2-2rr'\cos(d(x,x')\wedge\pi)},
            & K=0,
        \end{cases}
    \end{equation}
    and
    \begin{equation}
        \m_{\Sigma_K^N}(\d r,\d x)
        \coloneqq
        \sin_K^N(r)\,\m(\d x)\d r.
    \end{equation}
    We refer to $\Sigma_0(F)$ as the \textbf{Euclidean cone} over $F$ and to $\Sigma_K(F)$ with $K>0$ as the \textbf{spherical cone} over $F$.
    \end{definition}
Here the function $\sin_K$ is defined as follows
\[
\sin_K(r) \coloneqq 
\begin{cases}
	\frac{1}{\sqrt{K}} \sin(\sqrt{K}r) & \text{for } K > 0, \\
	r & \text{for } K = 0, \\
	\frac{1}{\sqrt{|K|}} \sinh(\sqrt{|K|}r) & \text{for } K < 0.
\end{cases}
\]
Let $(F,h,e^{-f}\vol_h)$ be an $n$-dimensional smooth weighted manifold. 
For $N\geq n$ and $K\in\R$, the regular part of its $(K,N)$-cone is a weighted Riemannian manifold with
\begin{equation}\label{eq:weighted_measure_cone}
	\m_{\Sigma^N_K}=e^{-\Phi(r,x)} \cdot \sin_K^n(r)\d \vol_h(x)\d r=e^{-\Phi(r,x)}\vol_{g_{\Sigma_K}}\,,
	\end{equation}
where 
\begin{equation}\label{eq:weight_of_cone}
	\Phi(r,x)\coloneqq f(x)-(N-n)\log \sin_K(r)\,.
\end{equation}
For each $(r,x)$ in the punctured cone and $V,W\in T_xF$, we denote their vertical lifts by $\widetilde V$ and $\widetilde W$, respectively, and identify $\partial_r$ with its horizontal lift. 
Thus,
\begin{equation}\label{eq:conic_g_sigma}
	g_{\Sigma_K}(\partial_r,\partial_r)=1,
	\qquad
	g_{\Sigma_K}(\partial_r,\widetilde V)=0,
	\qquad
	g_{\Sigma_K}(\widetilde V,\widetilde W)
	=
	\sin_K^2(r)h(V,W).
\end{equation}
Denote by $\Ric^{\Sigma}_{\Phi}$ the Bakry--\'Emery Ricci tensor on $(\Sigma_K(F),d_{\Sigma_K},\m_{\Sigma^N_{K}})$ and by $\Ric^F_f$ the Bakry--\'Emery Ricci tensor on $(F,h,e^{-f}\vol_h)$.
A direct computation using the warped-product formulas (cf. \cite[Proposition~3.2]{Ketterer-JFA13}) gives
\begin{align}
		\Ric^{\Sigma}_{\Phi}(\partial_r,\partial_r)&=nK+\frac{N-n}{\sin_K^2(r)}\,,\qquad 	\Ric^{\Sigma}_{\Phi}(\partial_r,\widetilde V)=-\frac{\sin_K'(r)}{\sin_K(r)}\,\d f(V)\,,\label{ineq:ONeill-1}\\
		\Ric_{\Phi}^{\Sigma}(\widetilde V,\widetilde W)&=	NK\sin_K^2(r)h(V,W)+
		\bigl(
		\Ric^F_f-(N-1)h
		\bigr)(V,W)\,.\label{ineq:ONeill-2}
\end{align}
Moreover, by the Bishop--O'Neill formula \cite{Bishop-ONeill69} (see also \cite[Proposition~1.5]{Ketterer26}), one has
\begin{align}\label{eq:Bishop-ONeill_Sec}
   \Sec^\Sigma(\partial_r,\widetilde V)=K,
  \qquad
  \Sec^\Sigma(\widetilde V,\widetilde W)
  =
  K+\frac{\Sec^F(V,W)-1}{\sin_K^2(r)}.
\end{align}

\section{Boundary rigidity via lower entropic Ricci bounds}\label{sec:boundary-rigidity}
The lower $\eta$-quantity $\eta_*(p)$ reflects the entropic behavior of the metric measure structure at arbitrarily small scales near $p$.
For $p\in\partial M$, the condition $\eta_*(p)>-\infty$ turns out to characterize precisely the local convexity of $\partial M$ near $p$.
The following result may therefore be viewed as a local version of Han's boundary rigidity theorem~\cite[Corollary~2.5]{Han-20}.

Denote by $\nu$ the inward unit normal vector field along $\partial M$.
The \textbf{second fundamental form} of $\partial M$ is
\begin{equation}
	\II_q(X,Y)\coloneqq-g_q(\nabla_X\nu,Y),\qquad q\in\partial M,\ X,Y\in T_q\partial M.
\end{equation}
We write $\II\geq0$ on $\Sigma\subset\partial M$ if $\II_q(X,X)\geq0$ for every $q\in\Sigma$ and $X\in T_q\partial M$.

\begin{theorem}[Boundary rigidity]\label{thm:boundary_rigidity}
	Let $(M,g,\m)$ be a smooth weighted Riemannian manifold with boundary, where $\m\coloneqq e^{-f}\vol_g$ with $f\in C^\infty(M)$.
	Let $p\in\partial M$.
	Then $\eta_*(p)>-\infty$ if and only if $\II\geq0$ on some neighborhood of $p$ in $\partial M$.
\end{theorem}

\subsubsection*{Collar neighborhood and extension across the boundary}
Fix a relative neighborhood $\Sigma$ of $p$ in $\partial M$.
After shrinking $\Sigma$ if necessary, there exists $R>0$ such that the geodesic collar
\begin{equation}\label{eq:normal_collar_map}
	\mathcal C\colon \Sigma\times[0,R)\to M,
	\qquad
	(q,r)\mapsto\exp_q(r\nu_q),
\end{equation}
is a diffeomorphism onto a neighborhood of $\Sigma$ in $M$; cf.~\cite{Schick01}.
By the Gauss lemma,
\begin{equation}\label{eq:collar_metric}
	{\mathcal C}^*g=\d r^2+h_r,
\end{equation}
where $(h_r)_{0\leq r<R}$ is a smooth family of Riemannian metrics on $\Sigma$; see~\cite[Proposition~2.8]{Schick01}.

Let $M^+\cup_{\partial M}M^-$ be the smooth double of $M$, with $M$ identified with $M^+$.
Via ${\mathcal C}$, we identify a neighborhood of $\Sigma$ in the double with $\Sigma\times(-R,R)$, where $r\geq0$ corresponds to $M^+$.
After possibly decreasing $R$, the family $h_r$ extends smoothly to $-R<r<R$ as a family of Riemannian metrics.
Hence $\widetilde g\coloneqq g$ on $M^\circ$ and $\widetilde g\coloneqq\d r^2+h_r$ on $\Sigma\times(-R,R)$ define a smooth Riemannian metric on the open set
\[
\widetilde M\coloneqq M^\circ\cup\bigl(\Sigma\times(-R,R)\bigr).
\]
Extending $f$ smoothly to $\widetilde f$ on $\widetilde M$, we obtain a weighted manifold $(\widetilde M,\widetilde g,\widetilde\m)$ without boundary, $\widetilde\m\coloneqq e^{-\widetilde f}\vol_{\widetilde g}$, which extends $(M,g,\m)$ across $\Sigma$.

With the above sign convention, $\nabla^{\widetilde g}r=\nu$ on $\Sigma$ and
\begin{equation}\label{eq:II-and-Hess}
	\Hess_{\widetilde g}r(X,Y)
	=-\II(X,Y)
	\qquad\text{on }\Sigma
\end{equation}
for all $X,Y\in T\partial M$.
The condition $\II\geq 0$ is often referred to as the convexity of the boundary in the following sense: if $\II\geq0$ on $\Sigma$, then every point of $\Sigma$ has a neighborhood $\Omega\Subset\widetilde M$ such that for all $x,y\in\Omega\cap M$ the minimizing $\widetilde g$-geodesic from $x$ to $y$ lies in $M$; see \cite[Theorem~1.2.1]{Wang_book14}.\medskip

The sufficiency part of \cref{thm:boundary_rigidity} follows from this local convexity and the consistency of $\eta$-quantities with lower Ricci bounds established in \cite{Erbar-Li-Schultz26}.
\begin{proof}[Proof of \cref{thm:boundary_rigidity}, $\Leftarrow$ part]
	Assume that $\II\geq 0$ on a relative open set $\Sigma$ in $\partial M$.
	Let $\widetilde M$ be the extension above, and denote by $\eta$ and $\widetilde{\eta}$ the corresponding quantities defined on $(M,g,\m)$ and $(\widetilde M,\widetilde g,\widetilde\m)$, respectively.
	As discussed above, there exists a neighborhood $\widetilde{\Omega}\Subset \widetilde{M}$ of $\Sigma$ s.t. for all $x,y\in \widetilde{\Omega}\cap M$ the minimizing $\tilde g$-geodesic is contained in $M$.
	As $\widetilde g$ coincides with $g$ on $M$, any $W_2$-geodesic in $\widetilde{\Omega}\cap M$ w.r.t. $g$ is also a $W_2$-geodesic in $\widetilde{\Omega}$ w.r.t. $\widetilde{g}$.
	Thus by definition
	\begin{equation}
		\eta(x,y)\geq \widetilde{\eta}(x,y),\qquad\forall x,y\in \widetilde{\Omega}\cap M,
	\end{equation}
	and by \eqref{eq:eta-smooth-consistency} the latter has a lower bound in $\widetilde \Omega$ in terms of $\Ric^{\widetilde M}_{\widetilde g}$.
	This implies that $\eta_*$ admits a finite lower bound on $\Sigma$.
\end{proof}

The necessity part relies on the following technical lemma, proved in \cref{app:boundary-transport}.
\begin{lemma}\label{lem:boundary-transport}
	Suppose that $p\in\partial M$ and that
	\begin{equation}\label{eq:negative-II}
		\II_p(v,v)<0
	\end{equation}
	for some unit vector $v\in T_p\partial M$.
	Then there are a relative neighborhood $\Sigma\subset \partial M$ of $p$, a smooth signed-distance function $\rho\colon\Sigma\to\R$ to a smooth hypersurface in $\partial M$ through $p$ whose unit normal is $v$ at $p$, and a smooth extension $(\widetilde M,\widetilde g)$ of $M$ across $\Sigma$ with the following properties.
	
	Set $v_q\coloneqq\nabla^{\partial M}\rho(q)$ for $q\in\Sigma$, where the sign of $\rho$ is chosen so that $v_p=v$, and $E_\pm(q,t)\coloneqq\exp^{\widetilde M}_q(\pm tv_q)$.
	For every sufficiently small $\delta>0$, set
	\begin{equation}\label{eq:xy}
		x\coloneqq\exp_p^{\partial M}(-\delta v/2),
		\qquad
		y\coloneqq\exp_p^{\partial M}(\delta v/2).
	\end{equation}
	There exist an open set $\Omega$ in $M$ containing the boundary geodesic from $x$ to $y$, neighborhoods $U_x,U_y\subset\Omega$ of $x,y$, and smooth functions $u_\pm$ on $\Omega\cap M^\circ$ such that:
	\begin{enumerate}[label=(\roman*)]
		\item\label{item:bdy_calculus1}
		there is $r_0>0$ such that the maps $E_\pm\colon\Sigma\times(0,r_0)\to M^\circ$ are diffeomorphisms onto their images, and $\Omega\cap M^\circ$ is contained in both images;
		
		\item\label{item:bdy_calculus2}
		on $\Omega\cap M^\circ$,
		\begin{equation}\label{eq:u-pm}
			u_\pm\bigl(E_\pm(q,t)\bigr)=\rho(q)\pm t,
		\end{equation}
		and $u_\pm$ extend continuously to $\Omega\cap\partial M$ with $u_\pm|_{\partial M}=\rho$.
		For each $q\in\Sigma$, the curve $t\mapsto E_\pm(q,t)$ is an integral curve of $\pm\nabla u_\pm$, i.e.
		\begin{equation}\label{eq:integral-curve-u_pm}
			\partial_tE_\pm(q,t)=\pm\nabla u_\pm\bigl(E_\pm(q,t)\bigr).
		\end{equation}
		Moreover, on $\Omega\cap M^\circ$,
		\begin{equation}\label{eq:laplacian-estimates}
			|\nabla u_\pm|=1, \qquad \Delta u_\pm\bigl(E_\pm(q,t)\bigr)=\pm\frac1t+O(1),
		\end{equation}
		where the $O(1)$-term is bounded uniformly in $(q,t)$;
		
		\item\label{item:bdy_calculus3}
		for every $\tilde x\in U_x\cap M^\circ$ and $\tilde y\in U_y\cap M^\circ$,
		\begin{equation}\label{eq:conjugate}
			d_g(\tilde x,\tilde y)\geq u_+(\tilde y)-u_-(\tilde x).
		\end{equation}
	\end{enumerate}
\end{lemma}

\textbf{Idea of the proof.}
When $\II$ is negative, minimizing geodesics between interior points can touch the boundary.
To illustrate this, assume $\II_p(v,v)<0$, and consider $\gamma_q$ for $q\in\partial M$ close to $p$ and $\delta>0$ sufficiently small that
\begin{equation}
	\gamma_q\colon [-\delta/2,\delta/2] \ni t\mapsto \exp_q^{\widetilde M}(tv_q).
\end{equation}
This curve touches $\partial M$ at $q$, lies in $M^\circ$ otherwise, and minimizes length between its endpoints. By \eqref{eq:integral-curve-u_pm}, its two halves follow the gradients of $u_-$ and $u_+$, respectively.

The proof turns this picture into constructing optimal transportations along a family of minimizing geodesics obtained by deforming $\gamma_q$, which travel from the interior to the boundary, follow the boundary, and return to the interior.
By \eqref{eq:conjugate} and \eqref{eq:u-pm}, the functions
\begin{equation}
	\varphi\coloneqq-\delta u_-,
	\qquad
	\psi\coloneqq\delta u_+-\delta^2/2
\end{equation}
form a pair of Kantorovich potentials for these transports\footnote{Notice that $\varphi$ and $\psi$ are conjugate w.r.t. $d^2_g/2$, i.e. $\varphi(x)+\psi(y)\leq d^2_g(x,y)/2$ for all $x,y$; and $\varphi(\tilde x)+\psi(\tilde y)=\frac12 d^2_g(\tilde x,\tilde y)$ for $\tilde x\coloneqq \gamma_q(-\delta/2)$ and $\tilde y\coloneqq \gamma_q(\delta/2)$}.
The initial and terminal transport velocities are therefore $-\nabla\varphi=\delta\nabla u_-$ and $\nabla\psi=\delta\nabla u_+$, respectively, and the endpoint derivatives of the entropy are governed by $\Delta u_-$ and $\Delta u_+$.
By the estimate \eqref{eq:laplacian-estimates}, these Laplacians are of order $\delta^{-1}$, producing an unbounded defect in entropy convexity.

\begin{proof}[Proof of \cref{thm:boundary_rigidity}, $\Rightarrow$ part]
	Suppose that \eqref{eq:negative-II} holds.
	Take the objects provided by \cref{lem:boundary-transport} and fix such a $\delta$.\medskip
	
	\textbf{Part 1:} construction of the transport.
	Let $\Phi_s$ be the flow of the unit vector field $q\mapsto v_q$ on $\Sigma$.
	Since $\rho$ is a signed-distance function, the flow lines are unit-speed geodesics of $\partial M$, i.e.\ $\Phi_s(q)=\exp^{\partial M}_q(sv_q)$, and
	\begin{equation}\label{eq:boundary-exp}
		\rho\bigl(\Phi_s(q)\bigr)=\rho(q)+s.
	\end{equation}
	In particular, $x=\Phi_{-\delta/2}(p)$ and $y=\Phi_{\delta}(x)$.
	
	Let $q\in\Sigma$ be close to $x$ and $\tau\in(0,r_0)$ be small.
	Then $q'\coloneqq\Phi_{\delta-2\tau}(q)\in\Sigma$, and we let $\gamma_{q,\tau}\colon[0,\delta]\to M$ be the curve (see \cref{fig:boundary-transport} below)
	\begin{equation}\label{eq:transport-geodesic}
		\gamma_{q,\tau}(s)\coloneqq
		\begin{cases}
			E_-(q,\tau-s), &0\leq s\leq\tau,\\
			\Phi_{s-\tau}(q), &\tau\leq s\leq\delta-\tau,\\
			E_+(q',s-\delta+\tau), &\delta-\tau\leq s\leq\delta.
		\end{cases}
	\end{equation}
	By \cref{lem:boundary-transport}~\ref{item:bdy_calculus1} and continuity, the first and the last pieces (except for $q$ and $q'$) lie in $U_x\cap M^\circ$ and $U_y\cap M^\circ$, respectively.
	The curve $\gamma_{q,\tau}$ is parametrized by arclength and has length $\delta$.
	Its endpoints are $z\coloneqq E_-(q,\tau)$ and
	\begin{equation}\label{eq:T}
		T(z)\coloneqq E_+(q',\tau).
	\end{equation}
	As a composition of $E_-^{-1}$, $(q,\tau)\mapsto(q',\tau)$, and $E_+$, the map $T$ is a diffeomorphism onto its image by \cref{lem:boundary-transport}~\ref{item:bdy_calculus1}.
	By \eqref{eq:u-pm} and \eqref{eq:boundary-exp},
	\begin{equation}\label{eq:potential-difference}
		u_+\bigl(T(z)\bigr)-u_-(z)=\rho(q')+\tau-\rho(q)+\tau=\delta.
	\end{equation}
	Hence $d_g(z,T(z))=\delta$ by \eqref{eq:conjugate}, i.e.\ $\gamma_{q,\tau}$ is a minimizing geodesic in $M$.
    \begin{figure}[htbp]
\centering
\begin{tikzpicture}[scale=2, line cap=round, line join=round, font=\small,
    >={Stealth[length=2mm]}]
    \definecolor{coneaccent}{RGB}{35,78,113}
    \def\R{2}
    \def\dl{4}
    \def\ta{0.7}
    \def\ep{0.8}
    \pgfmathsetmacro{\thx}{90+0.5*\dl/\R*180/pi}
    \pgfmathsetmacro{\thy}{90-0.5*\dl/\R*180/pi}
    \pgfmathsetmacro{\thq}{\thx-\ta/\R*180/pi}
    \pgfmathsetmacro{\thqq}{\thy+\ta/\R*180/pi}
    \begin{scope}[yscale=-1]
    \coordinate (C) at (0,-\R);
    \coordinate (p) at ($(C)+(90:\R)$);
    \coordinate (x) at ($(C)+(\thx:\R)$);
    \coordinate (y) at ($(C)+(\thy:\R)$);
    \coordinate (q) at ($(C)+(\thq:\R)$);
    \coordinate (qq) at ($(C)+(\thqq:\R)$);
    \coordinate (z) at ($(q)+(\thq+90:\ta)$);
    \coordinate (Tz) at ($(qq)+(\thqq-90:\ta)$);

    \fill[black!7] ($(C)+(172:\R)$) arc[start angle=172, end angle=8, radius=\R] -- cycle;
    \draw[line width=0.8pt] ($(C)+(172:\R)$) arc[start angle=172, end angle=8, radius=\R];
    \node at (2.2,0.45) {$M$};
    \node[black!65] at (0,-1.2) {$M^c$};
    \node[above right=-2pt] at ($(C)+(8:\R)$) {$\partial M$};

    \foreach \c in {x,y}
        \draw[black!55, densely dashed, line width=0.5pt] (\c) circle (\ep);

    \draw[coneaccent, line width=1.1pt, ->] (z) -- ($(z)!0.55!(q)$);
    \draw[coneaccent, line width=1.1pt] ($(z)!0.5!(q)$) -- (q);
    \draw[coneaccent, line width=1.1pt] (q) arc[start angle=\thq, end angle=\thqq, radius=\R];
    \draw[coneaccent, line width=1.1pt, ->] (qq) -- ($(qq)!0.55!(Tz)$);
    \draw[coneaccent, line width=1.1pt] ($(qq)!0.5!(Tz)$) -- (Tz);

    \draw[black!60, <->, line width=0.4pt]
        ($(C)+(\thq:\R+0.3)$) arc[start angle=\thq, end angle=\thqq, radius=\R+0.3];
    \node[below] at ($(C)+(90:\R+0.3)$) {$\delta-2\tau$};
    \node[coneaccent] at ($(C)+(101:\R-0.33)$) {$\gamma_{q,\tau}$};
    \node[below left=-1pt] at ($(z)!0.5!(q)$) {$\tau$};
    \node[below right=-1pt] at ($(qq)!0.5!(Tz)$) {$\tau$};

    \foreach \P in {p,x,y,q,qq} \fill (\P) circle (1.3pt);
    \fill[coneaccent] (z) circle (1.5pt);
    \fill[coneaccent] (Tz) circle (1.5pt);
    \node[above=2pt] at (p) {$p$};
    \node[above=2pt] at (q) {$q$};
    \node[above=2pt] at (qq) {$q'$};
    \node[above right=-1pt] at (x) {$x$};
    \node[above left=-1pt] at (y) {$y$};
    \node[below left=0pt] at (z) {$z$};
    \node[below right=1pt] at (Tz) {$T(z)$};
    \node[left] at ($(x)+(-\ep,-0.25)$) {$B_\varepsilon(x)$};
    \node[right] at ($(y)+(\ep,-0.25)$) {$B_\varepsilon(y)$};
    \end{scope}
    \node at (0,-0.85) {$z=E_-(q,\tau),\qquad T(z)=E_+(q',\tau)$};
\end{tikzpicture}
\vspace{-1pt}
\caption{Positions of the points chosen in Part~1.}
\label{fig:boundary-transport}
\end{figure}
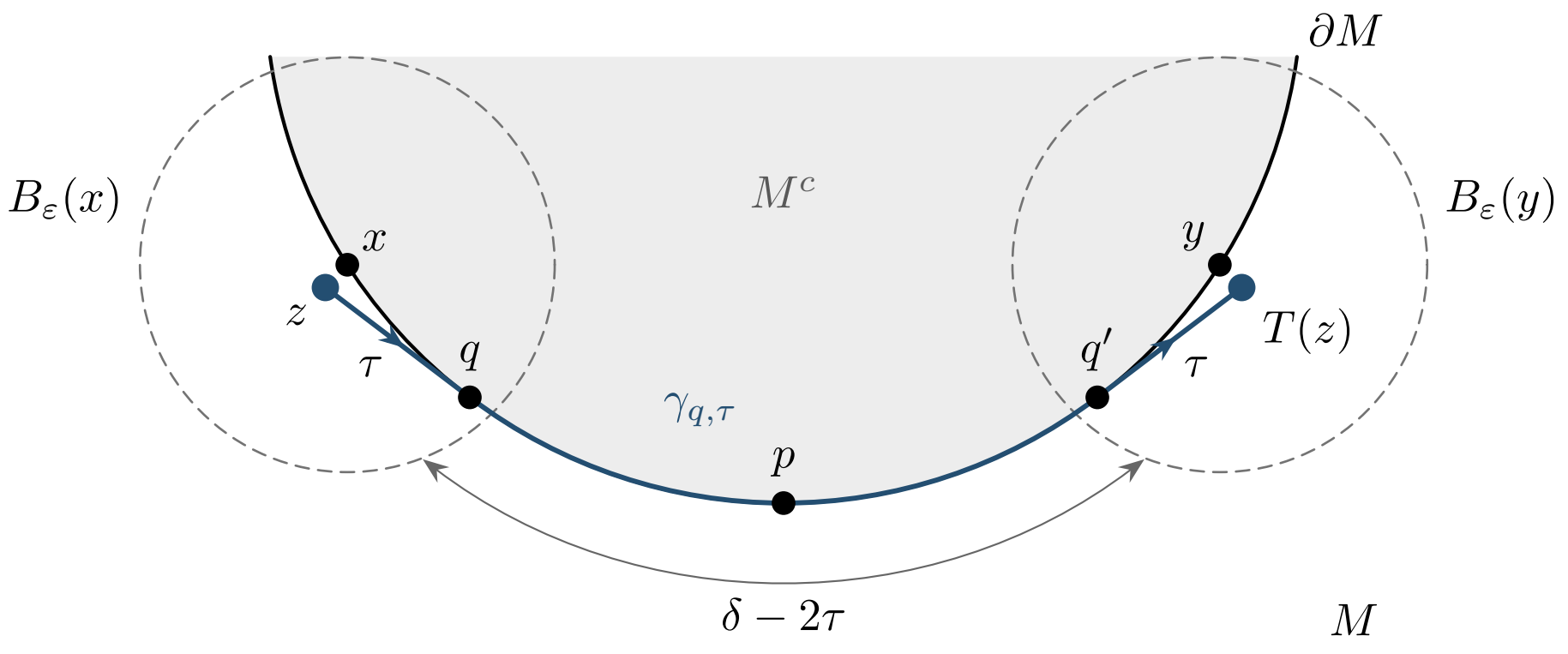
	Fix $\varepsilon>0$ and let $h>0$ be small.
	Since $T(z)\to y$ as $z\to x$, we can choose a probability measure $\mu^0$ with smooth density, supported in $\{E_-(q,\tau):h<\tau<2h\}\cap B_\varepsilon(x)$, s.t. $\mu^1\coloneqq T_\#\mu^0$ is supported in $B_\varepsilon(y)$.
	Both measures have smooth densities and hence finite entropy.
	Set $\varphi\coloneqq-\delta u_-$ and $\psi\coloneqq\delta u_+-\delta^2/2$.
	For $\tilde x\in\spt\mu^0$ and $\tilde y\in\spt\mu^1$, \eqref{eq:conjugate} and $(d_g-\delta)^2\geq0$ give
	\begin{equation}\label{eq:dual-calibration}
		\varphi(\tilde x)+\psi(\tilde y)\leq\delta\,d_g(\tilde x,\tilde y)-\frac{\delta^2}{2}\leq\frac12d_g^2(\tilde x,\tilde y),
	\end{equation}
	and by \eqref{eq:potential-difference} the equality holds if $\tilde y=T(\tilde x)$ and $d_g(\tilde x,T(\tilde x))=\delta$.
	Therefore, $(\varphi, \psi)$ is a pair of Kantorovich potentials for the optimal transport from $\mu^0$ to $\mu^1$, and $T$ is an optimal map.
	Consequently, $\mu^a\coloneqq(e_a)_\#\pi$ is a $W_2$-geodesic from $\mu^0$ to $\mu^1$, where $\pi$ is the push-forward of $\mu^0$ under $E_-(q,\tau)\mapsto\gamma_{q,\tau}(\delta\,\cdot)$.\medskip
	
	\textbf{Part 2: entropy derivative.}
	Since $\tau\in(h,2h)$ on $\spt\mu^0$, each curve $\gamma_{q,\tau}(\delta\,\cdot)$ runs on its first piece for $a\in[0,h/\delta]$ and on its last piece for $a\in[1-h/\delta,1]$ in $M^\circ$.
	By \eqref{eq:integral-curve-u_pm}, its velocity there is $\delta\nabla u_-$ and $\delta\nabla u_+$, respectively.
	Hence $(\mu^a)$ solves the continuity equation with these velocity fields on the two intervals, and the standard computation (see e.g.\ the proof of~\cite[Proposition~4.10]{Erbar-Li-Schultz26}) gives
	\begin{equation}\label{eq:entropy-slopes}
		\partial_a^+\ent(\mu^a)|_{a=0}=-\delta\int\Delta u_-\d\mu^0,
		\qquad
		\partial_a^-\ent(\mu^a)|_{a=1}=-\delta\int\Delta u_+\d\mu^1.
	\end{equation}
	The points of $\spt\mu^0$ are of the form $E_-(q,\tau)$ with $\tau\in(h,2h)$, and their images under $T$ are $E_+(q',\tau)$ with the same $\tau$.
	Thus \eqref{eq:laplacian-estimates} yields, with $C$ independent of $h$ and $\varepsilon$,
	\begin{equation}
		\partial_a^+\ent(\mu^a)|_{a=0}\geq\frac{\delta}{2h}-C\delta,
		\qquad
		\partial_a^-\ent(\mu^a)|_{a=1}\leq-\frac{\delta}{2h}+C\delta.
	\end{equation}
	For the weight, $\ent_\m(\mu)=\ent(\mu)+\int f\d\mu$.
	Since $\partial_a\int f\d\mu^a=\int\langle\nabla f(\gamma_a),\dot\gamma_a\rangle\d\pi(\gamma)$ and $|\dot\gamma_a|=\delta$, the one-sided derivatives of $a\mapsto\int f\d\mu^a$ at $a=0,1$ are bounded by $\delta\sup|\nabla f|$, where the supremum is taken over a compact neighborhood of the boundary geodesic from $x$ to $y$.
	Hence, for a possibly larger $C$,
	\begin{equation}
		\partial_a^-\ent_\m(\mu^a)|_{a=1}-\partial_a^+\ent_\m(\mu^a)|_{a=0}\leq-\frac{\delta}{h}+C\delta.
	\end{equation}
	Dividing by $W_2^2(\mu^0,\mu^1)=\delta^2$ and letting $h\to0$ yields $\eta_\varepsilon(x,y)=-\infty$.
	Since $\varepsilon>0$ is arbitrary, $\eta(x,y)=-\infty$.\medskip
	
	Suppose that $\eta_*(p)>-\infty$.
	By the definition of $\eta_*$ as a liminf, there are $r>0$ and $K\in\R$ such that $\eta(x',y')\geq K$ for all $x',y'\in B_r(p)$.
	If $\II_{p'}(v,v)<0$ for some $p'\in B_{r/2}(p)\cap\partial M$ and some unit vector $v\in T_{p'}\partial M$, then Parts~1 and~2, applied at $p'$ with $\delta<r$, produce $x,y\in B_r(p)$ with $\eta(x,y)=-\infty$, a contradiction.
	Therefore $\II\geq0$ on $B_{r/2}(p)\cap\partial M$.
\end{proof}

As a consequence of \cref{thm:boundary_rigidity}, we recover \cite[Corollary~2.5]{Han-20} via an approach based purely on optimal transport.
\begin{corollary}\label{cor:consistency_CD}
Let $(M,g)$ be a smooth complete Riemannian manifold with boundary and let $f$ be a smooth function on $M$. 
	Then for $K\in\R$ the following are equivalent 
	\begin{enumerate}
		\item\label{item:cor-CD1} $(M,d_g,e^{-f}\vol_g)$ satisfies $\rcd(K,\infty)$;
		\item\label{item:cor-CD2} $\eta_*|_{M^\circ}\geq K$, and $\eta_*|_{\partial M}>-\infty$;
        \item\label{item:cor-CD3} $\Ric_f\geq Kg$ on $M^\circ$, and $\II\geq 0$ on $\partial M$.
	\end{enumerate}
	\end{corollary}

\section{Upper entropic estimate in CAT spaces}\label{sec:CAT-upper}
As shown in the previous section, the presence of boundary may in general obstruct synthetic lower curvature bounds (a similar obstruction occurs for lower curvature bounds in the Alexandrov sense, since geodesics may branch at the boundary).
Nevertheless, every boundary point of a smooth Riemannian manifold with boundary admits a neighborhood with curvature bounded above in the Alexandrov sense; see \cite{ABB1993}.
In this section, we use such local upper curvature bounds to obtain upper estimates for the entropic quantity $\eta^*$.

Roughly speaking, a ${\rm CAT}(\kappa)$ space is a geodesic metric space in which geodesic triangles are not fatter than their comparison triangles in the model space of constant curvature $\kappa$.
We refer to, e.g. \cite[Chapter~II.1]{Bridson-Haefliger99} for the precise definition and basic properties.

Following the convention for $\sin_\kappa$ in \cref{sec:cones}, we denote
\[
\tan_\kappa(r)\coloneqq \frac{\sin_\kappa(r)}{\sin'_\kappa(r)}=\frac{1}{\sqrt{\kappa}}\tan(\sqrt{\kappa}r),\quad\text{for } \kappa>0;\qquad\tan_0(r)=r.
\]
\begin{theorem}\label{thm:CAT-eta-upper}
	Let $(M,g)$ be an $n$-dimensional Riemannian manifold, possibly with boundary.
	Assume that $(M,d)$ is a ${\rm CAT}(\kappa)$ space for some
	$\kappa\geq0$.
	Then for every $x\neq y\in M$ with $\ell\coloneqq d(x,y)<\pi/\sqrt{\kappa}$
	\begin{equation}\label{eq:eta-CAT}
		\eta(x,y)\leq n\kappa\frac{\tan_\kappa(\ell/2)}{\ell/2}.
	\end{equation}
	In particular, for every $p\in M$, $\eta^*(p)\leq n\kappa$.
\end{theorem}
 
\begin{proof}
	Let $x,y$ be distinct points in $M^\circ$ and let $\varepsilon>0$.
	Choose $0<\varepsilon_0<\min\{\varepsilon,\ell/4\}$ so small that, for $z=x,y$, we have $\overline{B_{2\varepsilon_0}(z)}\subset M^\circ$ and $\overline{B_{\varepsilon_0}(z)}$ lies in a normal neighborhood of $z$.
	If $\kappa>0$, also require $\ell+2\varepsilon_0<\pi/\sqrt\kappa$.\smallskip
	\\
   \textbf{Part 1:} endpoint measures and entropy derivatives.
	Choose $\chi\in C^\infty([0,\infty))$ with $0\leq\chi\leq1$, $\chi=1$ on $[0,\varepsilon_0/2]$, and $\chi=0$ on $[\varepsilon_0,\infty)$.
	For $z=x,y$ and $\delta\ll \varepsilon_0$, set
	\begin{equation}\label{eq:endpoint-measures}
		\rho_\delta^z(\cdot)\coloneqq \frac{\chi(d(z,\cdot))}{m_\delta^z}\exp\left(-\frac{d^2(z,\cdot)}{2\delta^2}\right)
		\qquad
		\mu_\delta^z\coloneqq\rho_\delta^z\vol_g,
	\end{equation}
	where $m_\delta^z$ normalizes the mass.
 
	Let $\mu^0_\delta \coloneqq \mu^x_\delta$ and $\mu^1_\delta\coloneqq \mu^y_\delta$.
	Let $\pi_\delta$ be the optimal dynamical plan from $\mu^0_\delta$ to $\mu^1_\delta$.
	Set $\mu_\delta^a\coloneqq(e_a)_\#\pi_\delta$ and $\sigma_\delta\coloneqq(e_0,e_1)_\#\pi_\delta$.
	Applying \cite[Theorem~23.14]{Villani} to short initial
	and reversed terminal segments gives
	\footnote{Although the global CD condition required in \cite[Theorem~23.14]{Villani} is not available here, \eqref{ineq:directional_derivative} still holds.
			Indeed, by the choice of the endpoint measures, there exists $h>0$ such that $\mu_\delta^a$ is supported in $\overline{B_{2\varepsilon_0}(x)}$ for all $a\in[0,h]$, and this ball lies in $M^\circ$, where the Ricci curvature is bounded from below.
			The same holds near $a=1$ by applying the theorem to the reversed transport.
	}
	\begin{align}
		\mathcal{DE}(\mu_\delta)
		&\coloneqq\partial_a^+\ent(\mu_\delta^a)|_{a=1-}
		-\partial_a^-\ent(\mu_\delta^a)|_{a=0+}\\
		&\leq\int\langle\nabla\log\rho_\delta^1(\gamma^1),\dot\gamma^1\rangle
		-\langle\nabla\log\rho_\delta^0(\gamma^0),\dot\gamma^0\rangle
		\d\pi_\delta(\gamma);\label{ineq:directional_derivative}
	\end{align}
	see also \cite[Lemma~8]{McCann-Topping}.
 
	Denote $D_{z}(\cdot)\coloneqq\frac12d^2(z,\cdot)$ for $z=x,y$.
	On $\{\rho_\delta^z>0\}$, we have
	\begin{equation}\label{eq:density-gradient}
		\nabla\log\rho_\delta^z=-\delta^{-2}\nabla D_{z}+R_\delta^z,
		\qquad
		R_\delta^z\coloneqq\nabla\log\chi(d(z,\cdot)).
	\end{equation}\smallskip
	\textbf{Part 2:} the leading term.
	For $v\in T_xM$ and $w\in T_yM$ with
	$|v|,|w|\leq\varepsilon_0$, set for $a\in [0,1]$
	\begin{equation}
		\mathtt{x}_v^a\coloneqq\exp_x(av),
		\quad
		\mathtt{y}_w^a\coloneqq\exp_y(aw),
		\quad
		L_{v,w}(a)\coloneqq
		\frac12d^2(\mathtt{x}_v^a,\mathtt{y}_w^a).
	\end{equation}
	Denote by $\gamma_{p,q}$ the minimizing geodesic from $p$ to $q$ parametrized on $[0,1]$.
	The first variation formula (we refer to \cite[Appendix]{ABB1993} for its proof on manifolds with boundary) gives
	\begin{align}\label{eq:first-variation}
		L_{v,w}'(0)
		&=\langle\dot\gamma_{x,y}(1),w\rangle-\langle\dot\gamma_{x,y}(0),v\rangle\\
		L_{v,w}'(1)
		&=\langle\nabla D_y(\mathtt{y}_w^1),\dot\gamma_{\mathtt{x}_v^1,\mathtt{y}_w^1}(1)\rangle-\langle\nabla D_x(\mathtt{x}_v^1),
		\dot\gamma_{\mathtt{x}_v^1,\mathtt{y}_w^1}(0)\rangle
	\end{align}
	where we used $\dot{\mathtt{x}}_v^1=\nabla D_x(\mathtt{x}_v^1)$ and
	$\dot{\mathtt{y}}_w^1=\nabla D_y(\mathtt{y}_w^1)$.
 
	Applying \cref{lemma:convexityCAT} to $L_{v,w}$, with $d(\mathtt{x}_v^a,\mathtt{y}_w^a)\leq \ell+2\varepsilon_0$, gives
	\begin{equation}\label{eq:CAT-cost-derivative-bound}
		L_{v,w}'(1)\geq L_{v,w}'(0)
		-\beta_\kappa(\ell+2\varepsilon_0)(|v|^2+|w|^2),
	\end{equation}
	where $\beta_\kappa(t)\coloneqq\sqrt\kappa\,t\tan\bigl(\sqrt\kappa\,t/2\bigr)$.
 
	Let $\widehat\sigma_\delta\coloneqq (\exp_x^{-1},\exp_y^{-1})_\#\sigma_\delta$ be the corresponding coupling on the tangent spaces. Then
	\begin{align}
		&-\int \langle \nabla D_y(\gamma^1),\dot\gamma^1\rangle -\langle \nabla D_x(\gamma^0),\dot\gamma^0\rangle \d\pi_\delta(\gamma)\\
		=&-\int \langle\nabla D_y(\mathtt{y}_w^1),\dot\gamma_{\mathtt{x}_v^1,\mathtt{y}_w^1}(1)\rangle-\langle\nabla D_x(\mathtt{x}_v^1),
		\dot\gamma_{\mathtt{x}_v^1,\mathtt{y}_w^1}(0)\rangle \d\widehat{\sigma}_\delta(v,w)\\
		=&-\int L_{v,w}'(1)\d\widehat{\sigma}_\delta(v,w)\leq  -\int \bigl(L_{v,w}'(0)- \beta_\kappa(\ell+2\varepsilon_0)(|v|^2+|w|^2)\bigr)\d\widehat{\sigma}_\delta(v,w)\\
		=& \int\langle\dot\gamma_{x,y}(0),v\rangle-\langle\dot\gamma_{x,y}(1),w\rangle\d\widehat{\sigma}_\delta(v,w)
		+\beta_\kappa(\ell+2\varepsilon_0)\int (|v|^2+|w|^2)\d\widehat{\sigma}_\delta(v,w)\label{eq:leading-term}
	\end{align}
	By \eqref{eq:moment-mu^z_delta} from \cref{lemma:moments-estimates} below,
	\begin{equation}\label{eq:moment-proof}
		\int |v|^2\d\widehat{\sigma}_\delta(v,w)=\int |w|^2\d\widehat{\sigma}_\delta(v,w)=n\delta^2+O(\delta^4).
	\end{equation}\smallskip
\textbf{Part 3:} error estimates.
	Since $\exp_z$ is a diffeomorphism around $B_{\varepsilon_0}(0)\subset T_zM$ for $z=x,y$, the measure $(\exp_z^{-1})_\#\vol_g$ there equals the volume measure induced by the pull-back metric $\exp_z^*g$ on $T_zM$.
	Thus, in normal coordinates,
	\begin{equation}\label{eq:normal-coordinate}
		(\exp_z^{-1})_\#\vol_g=\sqrt{\det{(g^z_{i,j}(v))}}\d \mathcal{L}^n(v),\qquad g^z_{i,j}(0)=\delta_{i,j},\quad \partial_kg^z_{i,j}(0)=0.
	\end{equation}
	Consequently, $(\exp_z^{-1})_\#\vol_g=(1+O(|v|^2))\mathcal L^n(v)$ uniformly in $\delta$.
 
	Since $\rho^0_\delta\circ\exp_x$ is radially symmetric and $\langle\dot\gamma_{x,y}(0),v\rangle$ is odd in $v$, we have
	\[
	\int \langle\dot\gamma_{x,y}(0),v\rangle\cdot \rho^0_\delta(\exp_x(v))\d\mathcal L^n(v)=0.
	\]
	Therefore, by \eqref{eq:m^z_delta} from \cref{lemma:moments-estimates} below,
	\begin{align}
		&\Bigl|\int \langle\dot\gamma_{x,y}(0),v\rangle\d\widehat{\sigma}_\delta(v,w)\Bigr|
		=\Bigl|\int \langle\dot\gamma_{x,y}(0),v\rangle\rho^0_\delta(\exp_x(v))\d (\exp_x^{-1})_\#\vol_g(v)\Bigr|
		\\
		\lesssim&
		\int |v|^3\rho^0_\delta(\exp_x(v))\d\mathcal L^n(v)\\
		\lesssim&
		\int (m_\delta^x)^{-1}e^{-|v|^2/(2\delta^2)}|v|^3\d\mathcal L^n(v)
		=O(\delta^3).\label{ineq:error-term1}
	\end{align}
	The same estimate holds also for $\int\langle\dot\gamma_{x,y}(1),w\rangle\d\widehat{\sigma}_\delta(v,w)$.
	Finally, for the remaining term $R^x_\delta$ (and similarly for $R^y_\delta$), we have by \eqref{eq:R^z_delta},
	\begin{equation}\label{ineq:R^x_delta}
		\Bigl|\int \langle R^x_\delta, \dot\gamma^0\rangle \d\pi_\delta(\gamma)\Bigr|\lesssim \int |R^x_\delta|\d\mu^0_\delta=o(1).
	\end{equation}\smallskip
	\textbf{Part 4:} conclusion.
	Combining the preceding parts, we conclude
	\begin{align}
		\mathcal{DE}(\mu_\delta)	&\overunderset{\eqref{ineq:directional_derivative}}{\eqref{eq:density-gradient}}{\leq}
		-\delta^{-2}\int\langle \nabla D_y(\gamma^1),\dot\gamma^1\rangle -\langle \nabla D_x(\gamma^0),\dot\gamma^0\rangle   \d\pi_\delta(\gamma) +
		\int  \langle R^y_\delta, \dot\gamma^1\rangle -\langle R^x_\delta, \dot\gamma^0\rangle \d\pi_\delta(\gamma)\\
		&\overunderset{\eqref{eq:leading-term}}{\eqref{ineq:R^x_delta}}{\leq}
		\delta^{-2}\beta_\kappa(\ell+2\varepsilon_0)\int (|v|^2+|w|^2)\d\widehat{\sigma}_\delta(v,w)\\
		&\qquad +\delta^{-2}\int\langle\dot\gamma_{x,y}(0),v\rangle-\langle\dot\gamma_{x,y}(1),w\rangle\d\widehat{\sigma}_\delta(v,w)+o(1)\\
		&\overunderset{\eqref{eq:moment-proof}}{\eqref{ineq:error-term1}}{=}
		2\delta^{-2}\beta_\kappa(\ell+2\varepsilon_0)(n\delta^2+O(\delta^4))+\delta^{-2}O(\delta^3)+o(1)\\
		&\,\,\,=\,\,\,2n\beta_\kappa(\ell+2\varepsilon_0)+o(1).
	\end{align}
	Since $\spt(\mu^z_\delta)\subset B_{\varepsilon}(z)$ for $z=x,y$ with $W_2(\mu_\delta^0,\mu_\delta^1)\xrightarrow{\delta\to 0} \ell$, we have by the definition of $\eta_\varepsilon$ that
	\begin{equation}
		\eta_\varepsilon(x,y)\leq
		\limsup_{\delta\to 0}
		\frac{\mathcal{DE}(\mu_\delta)}
		{W_2^2(\mu_\delta^0,\mu_\delta^1)}
		\leq
		\frac{2n\beta_\kappa(\ell+2\varepsilon)}{\ell^2}.
	\end{equation}
	Letting	$\varepsilon\downarrow0$ gives
	\begin{equation}
		\eta(x,y)\leq\frac{2n\beta_\kappa(\ell)}{\ell^2}=n\kappa\frac{\tan_\kappa(\ell/2)}{\ell/2}.
	\end{equation}
	The joint lower semi-continuity of $\eta$ on $M\times M$ extends this bound to all pairs of points.
	Finally, for any $p\in M$, letting $x,y\to p$ proves $\eta^*(p)\leq n\kappa$.
\end{proof}
 
\begin{lemma}\label{lemma:moments-estimates}
	Let $m^z_\delta$, $\mu^z_\delta$ and $R^z_\delta$ be given as in the proof of \cref{thm:CAT-eta-upper}, for $z=x,y$.
	Then as $\delta\to 0$
	\begin{align}
		m^z_\delta
		&=(2\pi)^{n/2}\delta^n+O(\delta^{n+2})\label{eq:m^z_delta}\\
		\int d^2(z,\cdot)\d\mu^z_\delta
		&=n\delta^2+O(\delta^4)\label{eq:moment-mu^z_delta}\\
		\int |R^z_\delta|\d\mu^z_\delta
		&= o(1).\label{eq:R^z_delta}
	\end{align}
\end{lemma}
\begin{proof}
	A standard computation shows
	\begin{equation}\label{eq:Gaussian-integrals}
		\int_{\R^n}e^{-|v|^2/(2\delta^2)}\d v=(2\pi)^{n/2}\delta^n,
		\qquad
		\int_{\R^n}|v|^2e^{-|v|^2/(2\delta^2)}\d v=n(2\pi)^{n/2}\delta^{n+2},
	\end{equation}
	and for any $N\in\N$
	\begin{equation}
		\int_{|v|\geq \varepsilon_0/2} e^{-|v|^2/(2\delta^2)}\d v=O(\delta^N).
	\end{equation}
	By the choice of the coordinate \eqref{eq:normal-coordinate}, there exist constants $c_z,C_z\in\R$ s.t.
	\[
	(1+c_z|v|^2)\mathcal L^n\leq (\exp_z^{-1})_\#\vol_g\leq (1+C_z|v|^2)\mathcal L^n,\qquad\text{on }B_{\varepsilon_0}(0)\subset T_zM.
	\]
	Since $m^z_\delta=\int\chi(|v|)\,e^{-|v|^2/(2\delta^2)}\d(\exp_z^{-1})_\#\vol_g(v)$, applying the above implies
	\[
	\int_{|v|\leq \varepsilon_0/2}e^{-|v|^2/(2\delta^2)}(1+c_z|v|^2)\d v\leq m^z_\delta\leq \int_{|v|\leq \varepsilon_0}e^{-|v|^2/(2\delta^2)}(1+C_z|v|^2)\d v,
	\]
	and hence \eqref{eq:m^z_delta} follows.
	Now, with \eqref{eq:m^z_delta} and \eqref{eq:Gaussian-integrals}, we have
	\begin{align}
		\int d^2(z,\cdot)\d\mu^z_\delta
		\leq& (m^{z}_\delta)^{-1} \int |v|^2 e^{-|v|^2/(2\delta^2)}(1+C_z|v|^2)\d v\\
		=&\bigl((2\pi)^{n/2}\delta^n+O(\delta^{n+2})\bigr)^{-1}\big(n(2\pi)^{n/2}\delta^{n+2}+O(\delta^{n+4})\big)=n\delta^2+O(\delta^4).
	\end{align}
	The reverse direction is similar and thus we obtain \eqref{eq:moment-mu^z_delta}.
 
	Finally, for \eqref{eq:R^z_delta}, as $\chi'\equiv 0$ on $[0,\varepsilon_0/2]$ and $|\nabla d(z,\cdot)|=1$, it follows by \eqref{eq:endpoint-measures}, \eqref{eq:density-gradient} and \eqref{eq:m^z_delta} that
	\begin{align}
		\int|R_\delta^z|\,\d\mu_\delta^z
		&=(m_\delta^z)^{-1}\int |\chi'(d(z,\cdot))|\, e^{-d^2(z,\cdot)/(2\delta^2)}\d\vol_g\\
		&\lesssim\delta^{-n} \int_{|v|\geq \varepsilon_0/2}e^{-|v|^2/(2\delta^2)}\d v=o(1).\tag*{\qedhere}
	\end{align}
\end{proof}
 
\begin{corollary}\label{cor:finite-at-boundary}
	Let $(M,g,e^{-f}\vol_g)$ be a weighted manifold, possibly with boundary.
	Then $\eta^*(p)<\infty$ for every $p\in M$.
\end{corollary}
\begin{proof}
	For $p\in M^\circ$, the claim follows from \cref{sec:eta_pre}.
	Let $p\in\partial M$ and assume first that $f\equiv0$.
	By \cite{ABB1993}, $(M,d)$ has curvature bounded above near $p$ in the sense of Alexandrov.
	Hence there are $\kappa\geq0$ and $s>0$ with $8s<\pi/\sqrt\kappa$ such that $B\coloneqq B_{4s}(p)$, endowed with the restriction of $d$, is a ${\rm CAT}(\kappa)$ space.
	In particular, $B$ is a Riemannian manifold with boundary whose intrinsic distance is $d$.
	For $x,y\in B_s(p)$ and $\varepsilon<s$, the triangle inequality shows that all minimizing geodesics between $B_\varepsilon(x)$ and $B_\varepsilon(y)$ lie in $B$, so $\eta(x,y)$ is the same for $M$ and for $B$.
	Thus \cref{thm:CAT-eta-upper} applied to $B$ gives $\eta^*(p)\leq n\kappa$.
 
	For general $f$, we have $\ent_\m(\mu)=\ent(\mu)+\int f\d\mu$, and near $a=0,1$ the transports in the proof of \cref{thm:CAT-eta-upper} run in $M^\circ$.
	Hence their entropy defect changes by $\int\langle\nabla f(\gamma^1),\dot\gamma^1\rangle-\langle\nabla f(\gamma^0),\dot\gamma^0\rangle\d\pi_\delta(\gamma)$, and it suffices to show that\footnote{In contrast to Part~2 of the proof of \cref{thm:boundary_rigidity}, where $\delta$ is fixed, the weight term has to be compared with $d^2(x,y)$ as $x,y\to p$; a bound by $\sup|\nabla f|$ alone only gives $O(d(x,y))$.}, for $x',y'\in M^\circ$ close to $p$ and every minimizing geodesic $\gamma\colon[0,1]\to M$ from $x'$ to $y'$,
	\begin{equation}\label{ineq:weight-term}
		(f\circ\gamma)'(1)-(f\circ\gamma)'(0)\leq C\,\ell^2,\qquad \ell\coloneqq d(x',y'),
	\end{equation}
	since then $\eta^*(p)\leq n\kappa+C$.
	Extend $f$ smoothly across $\partial M$, let $r$ be the collar coordinate of \cref{sec:boundary-rigidity}, i.e.\ the signed distance to $\partial M$, and set $F\coloneqq f+c\,r$ with $c>|\nabla f(p)|$, so that $\nabla F$ points strictly into $M$ near $p$.
	Let $C_F$ and $c_1$ be bounds for $|\Hess F|$ and $|\Hess r|$ near $p$.
 
	First, let $\Phi$ be the flow of $\nabla F$. For $\phi\in C^\infty_c((0,1))$ with $\phi\geq0$ and small $\lambda\geq0$, the curve $\gamma_\lambda\colon t\mapsto\Phi_{\lambda\phi(t)}(\gamma(t))$ stays in $M$ and has the same endpoints as $\gamma$.
	Since $\gamma$ minimizes the energy among such curves, differentiating at $\lambda=0$ (with $\partial_\lambda \gamma_\lambda|_{\lambda=0}=\phi\nabla F(\gamma)$ and $\nabla_\lambda\partial_t \gamma_\lambda=\nabla_t\partial_\lambda \gamma_\lambda$) gives
	\begin{align}
	0\leq \left.\frac{\d}{\d\lambda}\right|_{\lambda=0}\frac12\int |\dot\gamma_\lambda(t)|^2\d t&=\int^1_0 \langle \dot\gamma, \nabla_\lambda\partial_t\gamma_\lambda|_{\lambda=0}\rangle\d t  =\int^1_0 \langle \dot\gamma, \nabla_t \big( \phi\nabla F(\gamma) \big) \rangle\d t      \\
    &=\int_0^1\Big(\phi'\,(F\circ\gamma)'+\phi\,\Hess F(\dot\gamma,\dot\gamma)\Big)\d t,
	\end{align}
	that is, $(F\circ\gamma)''\leq C_F\ell^2$ in the sense of distributions on $(0,1)$.
 
	Second, $u\coloneqq r\circ\gamma\geq0$, and on $\{u>0\}$ the curve $\gamma$ is a geodesic of $g$, so 
    \[
    u''=\Hess r(\dot\gamma,\dot\gamma)+\langle \nabla r,\nabla_t\dot\gamma\rangle\geq-c_1|\dot\gamma|^2=-c_1\ell^2.
    \]
    As $f=F-c\,r$, we obtain $(f\circ\gamma)''\leq(C_F+c\,c_1)\ell^2$ in the sense of distributions.
	Since $\gamma$ is smooth near its endpoints, \eqref{ineq:weight-term} follows with $C\coloneqq C_F+c\,c_1$.
\end{proof}
\section{Entropic Ricci bounds in conic neighborhoods}
We now establish rigidity results for exact cones and conic neighborhoods using the Jacobian calculus and dilation arguments.
\subsection{Uniform curvewise coercivity bound under non-vanishing curvature}\label{sec:coercivity-Riem}
Let $(M,g)$ be a smooth $n$-dimensional Riemannian manifold, possibly with boundary.
\begin{definition}\label{def:geodesic-tube}
	Let $\Omega\subset M$ be open, and let $U_0,U_1\Subset\Omega$ be open.
	We call $\Omega$ a \textbf{geodesic tube} of $U_0$ and $U_1$ if every $(p,q)\in\overline U_0\times\overline U_1$ is joined by a unique minimizing geodesic $\sigma_{p,q}\colon[0,1]\to\Omega$, the union of their images is relatively compact
	in $\Omega$, and $\sigma_{p,q}^a$ depends smoothly on $(p,q,a)$.
\end{definition}

The following lemma will be convenient later. It follows from standard properties of points outside the cut locus together with a compactness argument, and it holds for possibly incomplete manifolds with boundary.
\begin{lemma}\label{lemma:geodesic-tube}
	Let $\gamma\colon [0,1]\to M^\circ$ be a minimizing geodesic from $p$ to $q$ s.t. $q\notin\Cut(p)$.
	Then there exist a neighborhood $V\subset T_pM$ of $\{a\dot\gamma(0):a\in[0,1]\}$ and neighborhoods $U_0,U_1\subset M$ of $p,q$, respectively, s.t. $\exp_p|_V$ is a diffeomorphism onto $\Omega\coloneqq\exp_p(V)$, and $\Omega$ is a geodesic tube of $U_0$ and $U_1$.
\end{lemma}

Let $U_0,U_1\Subset M$ be open and let $\Omega$ be a geodesic
tube of $U_0$ and $U_1$.
After shrinking $U_0$ if necessary, choose a smooth orthonormal frame $e_1(x),\ldots,e_n(x)$ on a neighborhood of $\overline U_0$.\footnote{This choice is only for convenience; the argument can equally be formulated using endomorphism fields along the geodesics.}

For $(x,y)\in U_0\times U_1$, let $E^i_{x,y}(a)$ be the parallel transport of $e_i(x)$ along $\sigma_{x,y}$, namely
\begin{equation}
	\nabla_a E^i_{x,y}=0,\qquad E^i_{x,y}(0)=e_i(x).
\end{equation}
In this frame, define the curvature matrix
\begin{equation}\label{eq:curvature-matrix}
	R_{x,y}(a)_{ij}\coloneqq
	\left\langle \Rm(E^j_{x,y}(a),\dot\sigma_{x,y}^a)\dot\sigma_{x,y}^a, E^i_{x,y}(a)\right\rangle.
\end{equation}
Both $E^i_{x,y}(a)$ and $R_{x,y}(a)$ depend smoothly on $(x,y,a)$.

\begin{lemma}[Uniform curvewise coercivity bound]\label{lemma:uniform_lower_bound}
	Let $q\notin \Cut(p)$ s.t. $p$ and $q$ are connected by a minimizing geodesic in $M^\circ$, along which $R_{p,q}$ is non-vanishing.
	Then there exist $\delta,c_0>0$, depending only on $R_{p,q}$, with the following property: if $R\colon[0,1]\to\R^{n\times n}$ is continuous with
		\begin{equation}\label{ineq:C0-close}
			\sup_{a\in[0,1]}|R(a)-R_{p,q}(a)|_{\HS}\leq\delta,
		\end{equation}
		then every absolutely continuous matrix-valued solution $S$ of
	\begin{equation}\label{eq:matrix-Riccati-R}
		S'+ S^2+ R=0
	\end{equation}
	on $[0,1]$ satisfies (independent of the initial value of $S$)
	\begin{equation}\label{ineq:coercivity-0}
		\int_0^1| S(a)|_{\HS}^2\d a\geq c_0.
	\end{equation}
	In particular, there exist neighborhoods $U_0, U_1$ of $p,q$ respectively s.t. for every $(x,y)\in U_0\times U_1$, \eqref{ineq:coercivity-0} holds for every absolutely continuous matrix-valued solution $S$ of
	\begin{equation}\label{eq:matrix-Riccati}
		S'+ S^2+ R_{x,y}=0
	\end{equation}
	along the geodesic $\sigma_{x,y}$.
\end{lemma}
\begin{proof}
	By \cref{lemma:geodesic-tube}, we can find neighborhoods $U_0$ and $U_1$ of $p$ and $q$ respectively and $\Omega\subset M^\circ$ s.t. $\Omega$ is a geodesic tube of $U_0$ and $U_1$.
	By assumption, there exists non-negative $\chi\in C_c^\infty((0,1))$ s.t.
	\begin{equation}
		\mathcal I\coloneqq\int_0^1\chi(a)| R_{p,q}(a)|_{\HS}^2\d a>0.
	\end{equation}
	Consider the fixed matrix field $B(a)\coloneqq\chi(a) R_{p,q}(a)$ along $\sigma_{p,q}$, and set
		\begin{equation}
			\alpha\coloneqq \| B'\|_{L^2((0,1);|{\cdot}|_{\HS})},
			\qquad
			\beta\coloneqq \|B\|_{L^\infty((0,1);|{\cdot}|_{\op})},
			\qquad
			\delta\coloneqq\frac{\mathcal I}{2\|B\|_{L^1((0,1);|{\cdot}|_{\HS})}}.
		\end{equation}
		If $R$ satisfies \eqref{ineq:C0-close}, then the Cauchy--Schwarz inequality gives
		\begin{equation}
			\int_0^1\langle R(a), B(a)\rangle_{\HS}\d a\geq\int_0^1\chi(a)| R_{p,q}(a)|_{\HS}^2\d a-\delta\|B\|_{L^1((0,1);|{\cdot}|_{\HS})}=\frac{\mathcal I}{2}.
	\end{equation}
	Now for any $ S$ as a solution to \eqref{eq:matrix-Riccati-R}, testing \eqref{eq:matrix-Riccati-R} against $B$, with integrating by parts and the Cauchy--Schwarz inequality, gives
	\begin{align}
		\frac{\mathcal I}{2}
		&\leq \int_0^1\langle R, B\rangle_{\HS}\d a= -\int^1_0 \langle S'+ S^2, B\rangle_{\HS} \d a\\
		&=\int_0^1\langle  S, B'\rangle_{\HS}\d a -\int_0^1\langle  S^2, B\rangle_{\HS}\d a\\
		&\le \alpha \left(\int^1_0 | S|^2_{\HS}\d a\right)^{1/2}+\beta\int^1_0 | S|^2_{\HS}\d a,\label{ineqeq:coercive}
	\end{align}
	where the boundary term vanishes since $B$ is compactly supported in $(0,1)$, and we have used
	\begin{equation}
		|\langle  S^2, B\rangle_{\HS}|
		\leq| B|_{\op}| S|_{\HS}^2.
	\end{equation}
	If $\| S\|_{L^2((0,1);|{\cdot}|_{\HS})}\geq1$, there is nothing to prove.
	If otherwise, \eqref{ineqeq:coercive} yields
	\begin{equation}
		\frac{\mathcal I}{2}\leq(\alpha+\beta)\| S\|_{L^2((0,1);|{\cdot}|_{\HS})}.
	\end{equation}
	As the constants $\mathcal I,\alpha,\beta,\delta$ depend only on $R_{p,q}$ and $\chi$, \eqref{ineq:coercivity-0} holds with $c_0\coloneqq\min\{1,\mathcal I^2/(4(\alpha+\beta)^2)\}$ uniform in $R$ and $ S$.
	Finally, since $R_{x,y}(a)$ depends smoothly on $(x,y,a)$, we can shrink $U_0,U_1$ s.t. $R_{x,y}$ satisfies \eqref{ineq:C0-close} for every $(x,y)\in U_0\times U_1$, which proves the last assertion.
\end{proof}

The Jacobian calculus turns this coercivity bound into positivity of $\eta$ when $\Ric_f\geq0$.
Note that all computations in \cref{sec:Jacobian-calculus} remain valid when restricted to a geodesic tube contained in the interior $M^\circ$.
\begin{proposition}\label{prop:positivity_eta}
	Let $(M,g,e^{-f}\vol_g)$ be a weighted manifold with $\Ric_f\geq 0$.
	Let $q\notin \Cut(p)$ s.t. $p$ and $q$ are connected by a minimizing geodesic in $M^\circ$, along which $R_{p,q}$ is non-vanishing.
	Then $\eta(p,q)>0$.
\end{proposition}
\begin{proof}
	Since $q\notin{\rm Cut}(p)$, one can find sufficiently small neighborhoods $U_0$ and $U_1$ of $p$ and $q$, respectively, admitting a geodesic tube.
	For any $W_2$-geodesic $(\mu^a)_{a\in[0,1]}$ with finite-entropy endpoints supported in $U_0,U_1$, respectively, \cref{lemma:Jacobiancalculus} applied within the tube gives
	\begin{align}
		\partial^{-}_a \ent_\m(\mu^a)\big|_{a=1-}-\partial_a^{+} \ent_\m(\mu^a)\big|_{a=0+}\geq  \int_{U_0}\int_0^1 | S_x(a)|_{\HS}^2 \d a\d\mu^0(x)\,,
	\end{align}
	where for $\mu^0$-a.e. $x$, $S_x$ solves the Riccati equation \eqref{eq:matrix-Riccati}.
	\cref{lemma:uniform_lower_bound} ensures that there exists a constant $c_0>0$ independent of the choices of $(\mu^a)$ s.t.
	\[
	\partial^{-}_a \ent_\m(\mu^a)\big|_{a=1-}-\partial_a^{+} \ent_\m(\mu^a)\big|_{a=0+}\geq c_0.
	\]
	Since $W_2(\mu^0,\mu^1)\leq d(p,q)+2\varepsilon$ whenever $\mu^0,\mu^1$ are supported in $B_\varepsilon(p)\subset U_0$ and $B_\varepsilon(q)\subset U_1$, the definition of $\eta$ yields $\eta(p,q)\geq c_0/d^2(p,q)>0$.
\end{proof}
\subsection{Rigidity of exact cones}

In this subsection, we establish a rigidity theorem for Euclidean cones in terms of the $\eta$-quantities.
Although the statement is restricted to $K=0$, this restriction is needed only for proving the rigidity of the sectional curvature of the fiber in the implication \ref{item:Econe-rigidity2}$\Rightarrow$\ref{item:Econe-rigidity3}, where the argument relies on the exact homogeneity of the Euclidean cone under dilations.
All other parts of the proof apply to general $K$.

\begin{theorem}[Rigidity of Euclidean cones]\label{thm:rigidity_E_cone}
	Let $(F,h,e^{-f}\vol_h)$ be an $n$-dimensional compact weighted manifold, possibly with boundary, where $n\geq2$.
	Let $\Sigma^N_K(F)$ be the $(K,N)$-cone of $F$ for $K=0$ and $N\geq n$, and let $o$ be the vertex.
	Then the following are equivalent:
	\begin{enumerate}[label=(\alph*)]
		\item\label{item:Econe-rigidity1} $\eta_*\equiv \eta^*\equiv nK$ on $\Sigma^N_K(F)$;
		\item\label{item:Econe-rigidity2} $\eta_*(o)>-\infty$, and $\eta^*(o)<\infty$;
		\item\label{item:Econe-rigidity3} $N=n$, $f$ is constant, $\Sec_h\equiv 1$, and $\II_{\partial F}\geq 0$.
	\end{enumerate}
\end{theorem}
\begin{proof}
	Write $\Sigma\coloneqq\Sigma_K^N(F)$, and let $g_\Sigma$ denote the Riemannian metric on its regular part.
	Let $\Ric^\Sigma_\Phi$ denote the Bakry--Émery Ricci tensor, with $\Phi$ given by \eqref{eq:weight_of_cone}.
	The implication \ref{item:Econe-rigidity1}$\Rightarrow$\ref{item:Econe-rigidity2} is trivial.
	
	\ref{item:Econe-rigidity2}$\Rightarrow$\ref{item:Econe-rigidity3}.
	By the definitions of $\eta_*$ and $\eta^*$, there exist $r_0>0$ and constants $K_-,K_+\in\R$ s.t. $K_-\leq\eta_*(p)\leq\eta^*(p)\leq K_+$ for all $p\in (0,r_0)\times F$.
	Then from \cref{sec:eta_pre}
	\begin{equation}
		K_-\,g_\Sigma\leq \Ric^\Sigma_\Phi\leq K_+\,g_\Sigma\qquad\text{on }(0,r_0)\times F.
	\end{equation}
	The radial and mixed components in \eqref{ineq:ONeill-1} force $N=n$ and $\d f\equiv0$, respectively. Since $F$ is connected, $f$ is constant.
	The tangential component in \eqref{ineq:ONeill-2} then gives $\Ric_h=(n-1)h$ on $F$.
	In particular, $\Ric_{g_\Sigma}=nK\,g_\Sigma$ on the punctured cone due to \eqref{eq:Bishop-ONeill_Sec}.
	
	We now show that $\Sec_h\equiv 1$ on $F$ under the assumption $K=0$.
	For $\lambda>0$, consider the dilation
	\begin{equation}
		\iota_\lambda\colon\Sigma\to\Sigma,\qquad \iota_\lambda(r,x)\coloneqq(\lambda r,x).
	\end{equation}
	Then
	\begin{equation}
		\iota_\lambda^*g_\Sigma=\lambda^2g_\Sigma,\qquad (\iota_\lambda)_\#\m_\Sigma=\lambda^{-(N+1)}\m_\Sigma.
	\end{equation}
	Consequently, for every $W_2$-geodesic $(\mu^a)_{a\in[0,1]}$ with finite entropy,
	\begin{equation}
		\ent((\iota_\lambda)_\#\mu^a)	=\ent(\mu^a)-(N+1)\log\lambda\,,
		\qquad
		W_2^2((\iota_\lambda)_\#\mu^0,(\iota_\lambda)_\#\mu^1)=\lambda^2W^2_2(\mu^0,\mu^1).
	\end{equation}
	Dilations bijectively map the admissible transports in the
	definition of $\eta_{\varepsilon/\lambda}(p,q)$ onto those in the definition of $\eta_\varepsilon(\iota_\lambda(p),\iota_\lambda(q))$, and leave the entropy derivatives unchanged.
	Thus $\eta_\varepsilon(\iota_\lambda(p),\iota_\lambda(q))=\lambda^{-2}\eta_{\varepsilon/\lambda}(p,q)$.
	Letting $\varepsilon\downarrow0$, we get
	\begin{equation}\label{eq:scaling_eta}
		\eta(\iota_\lambda(p),\iota_\lambda(q))=\lambda^{-2}\eta(p,q).
	\end{equation}
	Suppose that ${\Rm}^\Sigma\neq 0$ at $p=(r,x)$ for some $r>0$ and $x\in F^\circ$.
	Choose a unit vector $v\in T_p\Sigma$ with $\Rm^\Sigma(\cdot,v)v\neq0$, and let $q\coloneqq\exp_p(tv)$ for small $t>0$.
	Then $q\notin\Cut(p)$, the minimizing geodesic from $p$ to $q$ lies in the interior of the punctured cone, and along it the curvature matrix $R^\Sigma_{p,q}$ is not identically zero.
	Since the cone is Ricci flat, \cref{prop:positivity_eta} gives $\eta(p,q)>0$.
	Since $\iota_\lambda(p),\iota_\lambda(q)\to o$ as $\lambda\to0$, it follows from \eqref{eq:scaling_eta} that
	\begin{equation}
		\eta^*(o)\geq \limsup_{\lambda\to0}\eta(\iota_\lambda(p),\iota_\lambda(q))=\limsup_{\lambda\to0}\lambda^{-2}\eta(p,q)=+\infty.
	\end{equation}
	This contradicts \ref{item:Econe-rigidity2}.
	Hence $\Rm^\Sigma\equiv0$ on $(0,\infty)\times F^\circ$, and by continuity up to the boundary, the punctured cone is flat.
	The formula \eqref{eq:Bishop-ONeill_Sec} then forces $\Sec_h\equiv 1$ on $F$.
	
	It remains to prove that $\partial F$ is convex.
	Assume to the contrary that $\II_{\partial F}(w,w)<0$ for some $w\in T_x\partial F$, and fix $r\in(0,r_0)$.
	The inward unit normal to the cone boundary at $(r,x)$ is $\sin_K(r)^{-1}\widetilde\nu_F$, where $\widetilde\nu_F$ is the vertical lift of the inward unit normal vector $\nu_F$ to $\partial F$.
	Hence, for the vertical lift $\widetilde w\in T_{(r,x)}\Sigma$ of $w$,
	\begin{equation}
		\II_{\partial\Sigma}(\widetilde w,\widetilde w)=\sin_K(r)\,\II_{\partial F}(w,w)<0.
	\end{equation}
	By \cref{thm:boundary_rigidity}, $\eta_*((r,x))=-\infty$, which contradicts the assumption. \medskip
	
	\ref{item:Econe-rigidity3}$\Rightarrow$\ref{item:Econe-rigidity1}.
	By \ref{item:Econe-rigidity3} and \cref{cor:consistency_CD}, $(F,h,\vol_h)$ satisfies $\rcd(n-1,n)$.
	Then the $(K,n)$-cone $\Sigma$ satisfies $\rcd(Kn,n+1)$ (by \cite{Ketterer-JMPA15} and also \cite[Theorem~1.6]{Ketterer26-RCD-cone} for $K<0$), and consequently $\eta_*\geq nK$.
	
	It remains to prove $\eta^*\leq nK$.
	To this aim, we show that there exist $\ell_0,C>0$, depending only on $n$ and $K$, such that
	\begin{equation}\label{ineq:eta-upper-cone}
			\eta(p,q)\leq nK+C\,d_\Sigma^2(p,q)\qquad\text{whenever }0<d_\Sigma(p,q)<\ell_0.
		\end{equation}
	Since $\eta$ is jointly lower semi-continuous, it suffices to prove \eqref{ineq:eta-upper-cone} for a dense set of such pairs.
	As such pairs are dense, we may assume that $p,q\neq o$ are connected by a minimizing geodesic $(\gamma(a))_{a\in[0,1]}$ in the interior of the punctured cone, and that $q\notin\Cut(p)$.
	Denote $\ell\coloneqq d_{\Sigma}(p,q)$.
	By \cref{lemma:geodesic-tube}, there exists a neighborhood $V\subset T_p\Sigma$ of $\{a\dot\gamma(0):a\in[0,1]\}$ s.t. $\exp_p|_V $ is a diffeomorphism onto $\Omega\coloneqq \exp_p(V)$, and $\Omega$ is a geodesic tube of neighborhoods $U_p,U_q$ of $p,q$ respectively.
	
	Define, with $\exp_p^{-1}\coloneqq(\exp_p|_V)^{-1}$,
	\begin{equation}\label{eq:potential}
		\varphi(z)\coloneqq
		-g_p\bigl(\dot{\gamma}(0),\exp_p^{-1}z\bigr),\qquad z\in \Omega.
	\end{equation}
	Then $-\nabla\varphi(p)=\dot{\gamma}(0)$ and $\Hess\, \varphi(p)=0$.
	As the punctured cone has constant sectional curvature $K$, we have
	\begin{equation}
		\Hess\Big(\frac{d_\Sigma^2(\cdot, q)}{2}\Big)(p)
		=P^\parallel_\gamma+
		\ell\frac{\sin_K'(\ell)}{\sin_K(\ell)}P^\perp_\gamma,
	\end{equation}
	where $P^\parallel_\gamma$ and $P^\perp_\gamma$ stand for the projection of a vector $v\in T_p\Sigma$ to ${\rm span}(\dot\gamma(0))$ and its orthogonal complement; see e.g. \cite[Chapter~11]{Lee_Riem18}.
	Thus, choosing $\ell_0$ small and up to shrinking $U_p$ and $U_q$, we have
	\begin{equation}\label{eq:strict-cost-convepity}
		\Hess\Big(\frac{d_\Sigma^2(\cdot, y)}{2}-\varphi\Big)(x)\geq \frac12 {\rm Id}\qquad
		\text{for all }(x,y)\in U_p\times U_q.
	\end{equation}
	Let $T(z)\coloneqq\exp_z(-\nabla\varphi(z))$, so that $T(p)=q$, and shrink $U_p$ so that $T(U_p)\subset U_q$.
	Then following the argument of \cite[Theorem 13.5]{Villani}, $\varphi$ is $d^2_\Sigma$-concave relative to $U_p$ and $U_q$ (see \cite[Definition 3.1]{Cordero-Erausquin2001}), and $T$ is the optimal map from $\mu^0$ to $T_\#\mu^0$ for every probability measure $\mu^0\ll \vol_{g_\Sigma}|_{U_p}$.
	
	Now for the Jacobian calculus along transports induced by $T$, consider for any $x\in U_p$ the geodesic $\gamma_x\colon [0,1]\ni a\mapsto \exp_x(-a\nabla \varphi(x))$ and let $\ell_x\coloneqq d_\Sigma(x,T(x))$.
	Let $ A_x(a)$ be the Jacobi matrix along $\gamma_x$ given by \eqref{eq:Jacobi-matrix}.
	Then by \eqref{eq:Jacobi-equation} and the condition $\Sec_\Sigma\equiv K$ on $\Omega$,
	\begin{equation}
		\nabla_a \nabla_a  A_x+K\ell_x^2 P^\perp_{\gamma_x}\circ   A_x=0,\qquad  A_x(0)={\rm Id},\ \nabla_a A_x(0)=-\Hess \varphi(x).
	\end{equation}
	Direct computation yields, for $a\in[0,1]$,
	\[
	A_x(a)= P^\parallel_{\gamma_x} ({\rm Id}-a\cdot\Hess \varphi(x))+P^\perp_{\gamma_x}\big( \sin_K'(\ell_xa)\,{\rm Id}-\ell_x^{-1}\sin_K(\ell_x a)\, \Hess \varphi(x)\big),
	\]
	and in particular,
	\begin{equation}
		A_p(a)=P^\parallel_\gamma +\sin_K'(a\ell)P^\perp_\gamma.
	\end{equation}
	The explicit formula gives, uniformly for $a\in[0,1]$,
	\begin{equation}
		\|A_x(a)-{\rm Id}\|_{\rm op}+\|\nabla_a A_x(a)\|_{\rm op}\leq C\bigl(\ell_x^2+\|\Hess\varphi(x)\|_{\rm op}\bigr).
	\end{equation}
	Decreasing $\ell_0$ if necessary and after shrinking $U_p$, we have $\|A_x(a)^{-1}\|_{\rm op}\leq2$.
	Hence $S_x\coloneqq\nabla_a A_x\circ A_x^{-1}$, as in \eqref{eq:tensorS}, satisfies
	\begin{equation}\label{ineq:Hessian-estimate}
		\sup_{a\in[0,1]}\|S_x(a)\|_{\HS}
		\leq\sqrt{n+1}\sup_{a\in[0,1]}\|S_x(a)\|_{\rm op}
		\leq C\bigl(\ell_x^2+\|\Hess\varphi(x)\|_{\rm op}\bigr).
	\end{equation}
	Here and below, $C$ depends only on $n$ and $K$, independent of $(p,q)$ and chosen curve $(\mu^a)$.
	
	Take $\mu^0$ with a smooth compactly supported density in $U_p$, and let $(\mu^a)_{a\in[0,1]}$ be the induced $W_2$-geodesic, with $\mu^1=T_\#\mu^0$.
	Applying \cref{lemma:Jacobiancalculus} and \eqref{ineq:Hessian-estimate}, we obtain
	\begin{align}
		\partial^{-}_a \ent(\mu^a)|_{a=1-}-\partial_a^{+} \ent(\mu^a)|_{a=0+}
		&=\int_{U_p}\int_0^1
		\bigl(|S_x(a)|_{\HS}^2+nK|\dot\gamma_x(a)|^2\bigr)\d a\d\mu^0(x)\\
		&\leq nK W_2^2(\mu^0,\mu^1)+C\int_{U_p}
		\bigl(\ell_x^4+\|\Hess\varphi(x)\|^2_{\rm op}\bigr)\d\mu^0(x).
	\end{align}
	Shrinking $\mu^0$ to $\delta_p$, we have $\mu^1\to\delta_q$, $W_2^2(\mu^0,\mu^1)\to\ell^2$, and $\Hess\varphi(p)=0$.
	Consequently, 
	\begin{equation}
		\eta_\varepsilon(p,q)\leq nK+Cd_\Sigma^2(p,q),\qquad \forall \varepsilon>0
	\end{equation}
	    and letting $\varepsilon\downarrow0$ gives \eqref{ineq:eta-upper-cone}.
		By the reduction above, \eqref{ineq:eta-upper-cone} holds for all pairs with $0<d_\Sigma(p,q)<\ell_0$, including those involving the vertex or boundary points.
		Letting $p,q\to z$ gives $\eta^*(z)\leq nK$ for every $z\in\Sigma$.
	Together with $\eta_*\geq nK$, this proves \ref{item:Econe-rigidity1}.
\end{proof}

We give two remarks concerning the proof of \ref{item:Econe-rigidity3}$\Rightarrow$\ref{item:Econe-rigidity1}.
\begin{itemize}
	\item As in \cite[Section~4 and Lemma~A.3]{Erbar-Li-Schultz26}, the key point of the upper estimate for $\eta$ is to construct a Kantorovich potential with small Hessian.
	However, the construction in \cite{Erbar-Li-Schultz26} does not apply directly here, since the uniform neighborhood size in that construction depends on a positive lower bound for the injectivity radius (see \cite[Lemma~A.1]{Erbar-Li-Schultz26}), which is not available near the cone vertex.
	Here, we instead carry out the construction within a fixed geodesic tube, using that the sectional curvature is constant there.
	\item The upper estimate \cref{thm:CAT-eta-upper} does not apply here.
	Indeed, the round real projective space $F=\mathbb{RP}^n=\S^n/\mathbb{Z}_2$ has constant sectional curvature $1$ and empty boundary, but is not a ${\rm CAT}(1)$ space.
	Hence the Euclidean cone of $F$ is not a ${\rm CAT}(0)$ space by Berestovskii's theorem; see e.g. \cite[Theorem~II 3.14]{Bridson-Haefliger99}.
\end{itemize}
\subsection{Spaces with a conic singularity}\label{sec:conicnghd}
We give the precise definition of what we call a Riemannian manifold with a conic singularity.
\begin{definition}\label{def:conic-singularity}
	We say that $(M,d,g,o)$ is a Riemannian manifold with a \textbf{conic singularity} at $o$ \textbf{of order} $k\in\N$ if
	\begin{itemize}
		\item $(M,d)$ is a complete metric space;
		\item $(M_0,g)$ is an $n$-dimensional smooth Riemannian manifold where $M_0\coloneqq M\setminus\{o\}$,
		and the Riemannian distance induced by $g$ agrees with the restriction of $d$ to $M_0$;
		\item there exist a neighborhood $U_o$ of $o$, a number $l>0$, a compact smooth manifold $F$ \textbf{without boundary}, and a diffeomorphism
		\begin{equation}
			\Phi\colon (0,l)\times F\to U_o\setminus\{o\},\qquad \Phi^*g=\d r^2+r^2h_r
		\end{equation}
		where $(h_r)_{r\in(0,l)}$ is a smooth family of Riemannian metrics on $F$;
		\item there exist a smooth Riemannian metric $h_0$ on $F$ and some $\alpha>0$ such that, with
		\begin{equation}
			g_\Sigma\coloneqq \d r^2+r^2h_0,
		\end{equation}
		one has, uniformly on $F$ and for $j=0,1,\ldots, k$,
		\begin{equation}\label{eq:asymptotic_cone}
			\left|(\nabla^{g_\Sigma})^j \bigl(\Phi^*g-g_\Sigma\bigr)
			\right|_{g_\Sigma}=o(r^{\alpha-j}),\qquad r\to 0\,.
		\end{equation}
	\end{itemize}
\end{definition}
Here $k$ records the order of controlled derivatives, while $\alpha$ describes the decay rate.
The tensor norms in \eqref{eq:asymptotic_cone} are the pointwise norms induced by $g_\Sigma$; see, e.g. \cite[Proposition~2.40]{Lee_Riem18}.

For simplicity, we identify $U_o\setminus\{o\}$ with $(0,l)\times F$ via $\Phi$ and write $g$ for $\Phi^*g$.
We say that a family of tensor fields $(T_\lambda)_{\lambda>0}$ of a fixed type converges to $T_0$ in $C^k_{\mathrm{loc}}((0,\infty)\times F)$ as $\lambda\to0$ if, for every compact $Q\Subset(0,\infty)\times F$, $T_\lambda$ is defined on a neighborhood of $Q$ for all sufficiently small $\lambda$ and
\begin{equation}
	\sup_{Q}
	\left|(\nabla^{g_\Sigma})^j(T_\lambda-T_0)\right|_{g_\Sigma}
	\xrightarrow{\lambda\to0}0,
	\qquad 0\leq j\leq k.
\end{equation}
To study the geometry near $o$, we introduce the dilated metrics.
For $\lambda>0$, define
\begin{equation}
	\iota_\lambda\colon (0,l/\lambda)\times F\to (0,l)\times F,
	\qquad \iota_\lambda(r,x)\coloneqq(\lambda r,x),
\end{equation}
and set
\begin{equation}\label{eq:g_lambda}
	g_\lambda\coloneqq\lambda^{-2}\iota_\lambda^*g.
\end{equation}

\begin{lemma}[Convergence of dilated metrics]\label{lemma:Convergence_dilated_metrics}
	Let $(M,d,g,o)$ be an $n$-dimensional Riemannian manifold with a conic singularity at $o$ of order $k$, where $n\geq3$.
	Then $g_\lambda\to g_\Sigma$ in $C^k_{\mathrm{loc}}((0,\infty)\times F)$.
\end{lemma}
\begin{proof}
	Since the dilation changes only the radial coordinate and leaves the $F$-coordinate unchanged,
	\[
	\iota^*_\lambda g=\lambda^2 \d r^2+ \lambda^2 r^2 h_{\lambda r}\,.
	\]
	Thus
	\begin{equation}\label{eq:g_lambda-2}
		g_\lambda=\lambda^{-2}\iota_\lambda^*g=\d r^2+r^2h_{\lambda r}.
	\end{equation}
	As $g_\Sigma$ is conical, $\iota_\lambda^*g_\Sigma=\lambda^2g_\Sigma$.
	Therefore, putting $H\coloneqq g-g_\Sigma$, we have
	\begin{equation}\label{eq:Covar-deri-H}
		g_\lambda-g_\Sigma=\lambda^{-2}\iota_\lambda^*g-\lambda^{-2}\iota_\lambda^*g_\Sigma=\lambda^{-2}\iota_\lambda^*H.
	\end{equation}
	Under pullback, the Levi--Civita connection satisfies $\iota_\lambda^*(\nabla^{g_\Sigma}H)=\nabla^{\iota_\lambda^*g_\Sigma}(\iota_\lambda^* H)$.
	Moreover, a constant scaling of a metric does not change its Levi--Civita connection.
	Thus, with \eqref{eq:Covar-deri-H},
	\begin{equation}
		(\nabla^{g_\Sigma})^j(g_\lambda-g_\Sigma)= \lambda^{-2}(\nabla^{g_\Sigma})^j (\iota_\lambda^*H)=
		\lambda^{-2}\iota_\lambda^*
		\bigl((\nabla^{g_\Sigma})^jH\bigr),
		\qquad 0\leq j\leq k.
	\end{equation}
	Since $(\nabla^{g_\Sigma})^jH$ is a covariant tensor of rank $j+2$ and $\iota_\lambda^*g_\Sigma=\lambda^2g_\Sigma$,
	\begin{equation}\label{eq:scaled-cone-derivatives}
		\left|
		(\nabla^{g_\Sigma})^j(g_\lambda-g_\Sigma)
		\right|_{g_\Sigma}(s,z)
		=
		\lambda^j
		\left|
		(\nabla^{g_\Sigma})^jH
		\right|_{g_\Sigma}(\lambda s,z).
	\end{equation}
	Fix a compact annulus $A_{a,b}\coloneqq[a,b]\times F\subset(0,\infty)\times F$.
	For $0\leq j\leq k$, by \eqref{eq:asymptotic_cone} and \eqref{eq:scaled-cone-derivatives},
	\begin{align}
		\sup_{A_{a,b}}\left|
		(\nabla^{g_\Sigma})^j(g_\lambda-g_\Sigma)
		\right|_{g_\Sigma}
		&=\sup_{(r,z)\in A_{a,b}}\lambda^j
		\left|
		(\nabla^{g_\Sigma})^jH
		\right|_{g_\Sigma}(\lambda r,z)=o(\lambda^\alpha).
	\end{align}
	This completes the proof.
\end{proof}

\begin{proposition}\label{prop:lower_Ricci_cone}
	\label{prop:equiv_lower_bounds_conic}
	Let $(M,d,g,o)$ be an $n$-dimensional Riemannian manifold with a conic singularity at $o$ of order $0$, where $n\geq3$.
	Then 
	\[
	\eta_*(o)>-\infty \qquad \Leftrightarrow \qquad \Ric_g\geq Kg \text{ on $(0,l_0)\times F$ for some $K\in\R$ and $l_0>0$}.
	\]
	If in addition the singularity at $o$ is of order $2$, then these conditions imply $\Ric_{h_0}\geq(n-2)h_0$.
\end{proposition}

\begin{proof}
	Assume first that $\eta_*(o)>-\infty$ and fix $K<\eta_*(o)$.
	By definition, there exists $l_0>0$ s.t. $\eta(p,p')\geq K$	for all distinct $p,p'\in(0,l_0)\times F$.
	From \cref{sec:eta_pre}, $\Ric_g\geq Kg$ on $(0,l_0)\times F$.
	
	The reverse direction follows from \cite{Bertrand-Ketterer-Mondello-Richard21}.
	Indeed, shrink the conic neighborhood, leave $g$ unchanged near $o$, and modify it smoothly in an outer annulus so that it becomes a product near the outer boundary.
	Then doubling across this boundary produces a compact stratified space $(\widehat M,\widehat g)$ in the sense of \cite[Definition~1.1]{Bertrand-Ketterer-Mondello-Richard21} s.t.
	\begin{itemize}
		\item the regular part is a smooth manifold on which $\Ric_{\widehat g}\geq \widehat K \widehat g$ for some $\widehat K\in\R$;
		\item the two vertices, as the stratum of dimension $0$, form the singular part, which has codimension $n\geq 3$;
		\item the Riemannian metric $\widehat g$ on the regular part is an iterated edge metric in the sense of \cite[Definition~1.4]{Bertrand-Ketterer-Mondello-Richard21}, due to \eqref{eq:asymptotic_cone}.
	\end{itemize}
	Therefore, \cite[Theorem A]{Bertrand-Ketterer-Mondello-Richard21} gives that $\widehat M$, endowed with the induced distance $d_{\widehat g}$ and $\vol_{\widehat g}$, is an $\rcd(\widehat K, n)$ space.
	As $\widehat g$ coincides with $g$ around $o$ and $\eta$ is local, this implies that $\eta_*(o)\geq\widehat K>-\infty$.
	
	Now assume that the singularity is of order $2$ and $\Ric_g\geq Kg$ on $(0,l_0)\times F$.
	Consider the family of metrics $(g_\lambda)_{\lambda>0}$ given by \eqref{eq:g_lambda}.
	Since a constant rescaling of a metric does not change its Ricci tensor, we have
	\begin{equation}\label{ineq:Ric_g_lambda-1}
		\Ric_{g_\lambda}= \Ric_{\lambda^{-2}\iota^*_\lambda g}=\Ric_{\iota^*_\lambda g}=\iota_\lambda^*\Ric_g.
	\end{equation}
	Then with \eqref{ineq:Ric_g_lambda-1} and \eqref{eq:g_lambda}, we have
	\begin{equation}\label{ineq:Ric_g_lambda-2}
		\Ric_{g_\lambda}=\iota_\lambda^*\Ric_g\geq K \iota_\lambda^*g=K\lambda^2g_\lambda\qquad\text{on }(0,l_0/\lambda)\times F.
	\end{equation}
	By \cref{lemma:Convergence_dilated_metrics}, $g_\lambda\to g_\Sigma$ in $C^2_{\mathrm{loc}}((0,\infty)\times F)$, which yields that the curvature tensors converge in $C^0_{\mathrm{loc}}$.
	Passing to the limit in \eqref{ineq:Ric_g_lambda-2} gives $\Ric_{g_\Sigma}\geq0$.
	By the O'Neill formula \eqref{eq:Bishop-ONeill_Sec}, $\Ric_{h_0}\geq(n-2)h_0$.
	This completes the proof.
\end{proof}

We say that the conic singularity at $o$ is \textbf{removable} if $o$ is a topological manifold point of $M$, namely, if it admits a neighborhood homeomorphic to an open subset of
$\mathbb R^n$.

\begin{lemma}\label{lemma:removable_singularity}
	Let $(M,d,g,o)$ be an $n$-dimensional Riemannian manifold with a conic singularity at $o$ of order $0$, where $n\geq3$.
	Assume $\Sec_{h_0}\equiv 1$.
	Then the following are equivalent:
	\begin{enumerate}
		\item\label{item:remove_sing1} the singularity at $o$ is removable;
		\item\label{item:remove_sing2} $(F,h_0)$ is isometric to the unit sphere $\mathbb{S}^{n-1}$ with the round metric;
		\item\label{item:remove_sing3} the volume density at $o$ equals one:
		\begin{equation}
			\Theta(o)\coloneqq \lim_{r\to 0}\frac{\vol_g(B_r(o))}{\omega_n r^n}=1.
		\end{equation}
	\end{enumerate}
\end{lemma}
\begin{proof}
	Suppose first that $o$ is removable.
	Let $W\subset U_o$ be a neighborhood of $o$ homeomorphic to an open Euclidean ball.
	Then $W\setminus\{o\}$, being homeomorphic to a punctured ball, is simply connected since $n\geq 3$.
	We claim that $F$ is also simply connected.
	Indeed, fix $a>0$ with $B_{a}(o)\subset W$, and let $\gamma$ be a loop in $F$.
	Regarded as a loop in the slice $\{a/2\}\times F\subset W\setminus\{o\}$, it is null-homotopic there.
	Then composing this homotopy with the projection $(0,l)\times F\to F$ shows that $\gamma$ is null-homotopic in $F$.
	Since $F$ is compact and simply connected without boundary, applying the Killing--Hopf theorem (see e.g. \cite[Corollary 12.5]{Lee_Riem18}) with the property $\Sec_{h_0}\equiv 1$ gives that $(F,h_0)$ is isometric to the round sphere $\mathbb{S}^{n-1}$.
	This shows $\eqref{item:remove_sing1}\Rightarrow \eqref{item:remove_sing2}$.
	
	Assume that $(F,h_0)$ is the round sphere $\mathbb{S}^{n-1}$.
	Then the limiting cone with $g_\Sigma$ is the Euclidean space $\R^{n}$, and the polar identification $\Psi\colon U_o\to\R^n$, $\Psi(o)\coloneqq0$, $\Psi(r,x)\coloneqq rx$, is a homeomorphism onto the open ball of radius $l$.
	Moreover, by the zeroth-order estimate in \eqref{eq:asymptotic_cone}, $\Psi$ is bi-Lipschitz near $o$, and $g$ extends continuously across $o$ by $g|_o\coloneqq g_{\R^n}$.
	This shows $\eqref{item:remove_sing2}\Rightarrow \eqref{item:remove_sing1}$.
	
	We next compute the volume density for general $(F,h_0)$.
	Since $g-g_\Sigma=r^2(h_r-h_0)$ is purely tangential, the zeroth-order estimate in \eqref{eq:asymptotic_cone} gives
	\begin{equation}
		o(r^\alpha)=	\sup_{x\in F}|g-g_\Sigma|_{g_\Sigma}(r,x)=\sup_{x\in F} |h_r-h_0|_{h_0}(x).
	\end{equation}
	Therefore, uniformly on $F$,
	\begin{equation}
		\d\vol_{h_r}	=	\bigl(1+o(r^\alpha)\bigr)\d\vol_{h_0},
	\end{equation}
	and, since $B_\rho(o)\setminus\{o\}=(0,\rho)\times F$ for $\rho<l$,
	\begin{align}
		\vol_g(B_\rho(o))&=	\int_0^\rho r^{n-1}\vol_{h_r}(F)\d r =\frac{\vol_{h_0}(F)}{n}\rho^n+o(\rho^{n+\alpha}).
	\end{align}
	Hence $\Theta(o)=\vol_{h_0}(F)/\vol(\mathbb S^{n-1})$, where $\vol(\mathbb S^{n-1})=n\omega_n$. 
	Again by the Killing--Hopf theorem (see e.g. \cite[Theorem 12.4]{Lee_Riem18}), $	(F,h_0)\simeq\mathbb S^{n-1}/\Gamma$ for a finite group $\Gamma\leq O(n)$ acting freely on $\mathbb S^{n-1}$.
	By the computation above, $\Theta(o)=\vol(\mathbb S^{n-1}/\Gamma)/\vol(\mathbb S^{n-1})=1/|\Gamma|$.
	Therefore $\Theta(o)=1$ if and only if $\Gamma=\{\mathrm{id}\}$, i.e. $(F,h_0)\simeq\mathbb S^{n-1}$.
\end{proof}
\subsection{Proof of the main theorems}\label{sec:rigidity-conic-nghd}
\begin{proof}[Proof of \cref{thm:intro-conic}]
	Let $F$, $(h_r)_{r\in(0,l)}$ and $h_0$ be as in \cref{def:conic-singularity}, and let $g_\lambda=\d r^2+r^2h_{\lambda r}$ be the dilated metrics given by \eqref{eq:g_lambda} and \eqref{eq:g_lambda-2}.
	By \cref{lemma:Convergence_dilated_metrics}, $g_\lambda\to g_\Sigma$ in $C^2_{\mathrm{loc}}((0,\infty)\times F)$ as $\lambda\to 0$.
	In particular, the Christoffel symbols converge in $C^1_{\mathrm{loc}}$ and the curvature tensors converge in $C^0_{\mathrm{loc}}$.
	Moreover, by \cref{prop:lower_Ricci_cone}, $\Ric_g\geq Kg$ on $(0,l_0)\times F$ for some $l_0>0$ and $K\in\R$; hence \eqref{ineq:Ric_g_lambda-2} gives $\Ric_{g_\lambda}\geq K\lambda^2g_\lambda$ on $(0,l_0/\lambda)\times F$.
	
	Suppose that $\Sec_{h_0}\not\equiv 1$.
	Then by \eqref{eq:Bishop-ONeill_Sec}, $\Rm_{g_\Sigma}\neq0$ at some point of $(0,\infty)\times F$.
	As in the proof of \cref{thm:rigidity_E_cone}, we find a sufficiently short minimizing $g_\Sigma$-geodesic
	\begin{equation}
		\bar\gamma\colon[0,1]\to(0,\infty)\times F
	\end{equation}
	with endpoints $p,q$ s.t. $q\notin\Cut(p)$ and its curvature matrix $\overline{R}$, defined as in \eqref{eq:curvature-matrix}, is not identically zero.
	By \cref{lemma:geodesic-tube}, we may choose relatively compact neighborhoods $U_0,U_1$ of $p,q$, respectively, and a relatively compact open set $\Omega$ s.t.
	\begin{equation}
		\overline{U_0}\cup\overline{U_1}\Subset\Omega\Subset(0,\infty)\times F
	\end{equation}
	and $\Omega$ is a $g_\Sigma$-geodesic tube of $U_0$ and $U_1$.
	
	The contradiction argument is completed in the following two steps.\smallskip
	\\
	\textbf{Step 1}: Uniform coercivity along the dilated geodesic tubes.
	
	Due to the $C^2$-convergence $g_\lambda\to g_\Sigma$, the geodesic sprays (i.e. the vector fields on the tangent bundle whose integral curves encode the geodesics of the metric) of $g_\lambda$ converge in $C^1$ to that of $g_\Sigma$.
	Thus, after reducing $U_0,U_1$ and taking $\lambda$ sufficiently small, $\Omega$ is also a $g_\lambda$-geodesic tube of $U_0$ and $U_1$.
	In particular, every $(x,y)\in U_0\times U_1$ is joined by a unique minimizing $g_\lambda$-geodesic $\gamma_{\lambda,x,y}$ contained in $\Omega$, which depends continuously on $(\lambda,x,y)$.
	
	Put $g_0\coloneqq g_\Sigma$.
	For each $\lambda\geq0$, choose a $g_\lambda$-orthonormal frame on $U_0$ depending continuously on $\lambda$, and parallel transport it along $\gamma_{\lambda,x,y}$.
	Denote by $ R_{\lambda,x,y}$ the curvature matrix w.r.t. the corresponding metric and frame.
	The $C^1$-convergence of the Christoffel symbols and the $C^0$-convergence of the curvature tensors imply that $R_{\lambda,x,y}\to R_{0,x,y}$ uniformly in $(x,y,a)\in U_0\times U_1\times[0,1]$ as $\lambda\to0$.
	
	Since $R_{0,p,q}=\overline R\not\equiv0$, \cref{lemma:uniform_lower_bound} applied to $g_\Sigma$ yields constants $\delta,c_0>0$, depending only on $R_{0,p,q}$, s.t. \eqref{ineq:coercivity-0} holds for every absolutely continuous solution of \eqref{eq:matrix-Riccati-R} whenever $R$ satisfies \eqref{ineq:C0-close} with $R_{p,q}$ replaced by $R_{0,p,q}$.
		By the uniform convergence above and the continuity of $R_{0,x,y}$ in $(x,y)$, after shrinking $U_0,U_1$ and choosing $\lambda_0>0$ small, we have
		\[
		\sup_{a\in[0,1]}|R_{\lambda,x,y}(a)-R_{0,p,q}(a)|_{\HS}\leq\delta
		\qquad\text{for all }\lambda<\lambda_0,\ (x,y)\in U_0\times U_1.
		\]
		Note that this uses only the $C^0$-convergence of the curvature matrices.
		Hence every absolutely continuous solution $S$ of
	\begin{equation}\label{eq:Riccati_proof1}
		S'+ S^2+ R_{\lambda,x,y}=0
	\end{equation}
	along $\gamma_{\lambda,x,y}$ satisfies, independently of $\lambda$, $(x,y)$ and the initial value of $S$,
	\begin{equation}\label{ineq:uniform-Riccati-coercivity}
		\int^1_0 | S(a)|^2_{\rm HS}\d a\geq c_0.
	\end{equation}
	By the same argument as in \cref{prop:positivity_eta}, for any $W_{2,g_\lambda}$-geodesic $(\mu^a_\lambda)_{a\in[0,1]}$ whose endpoints have finite entropy and are supported in $U_0$ and $U_1$, respectively,
	\begin{align}
		&\mathcal{DE}_{\vol_{g_\lambda}}(\mu_{\lambda})\coloneqq\partial^{-}_a \ent_{\vol_{g_\lambda}}(\mu^a_\lambda)\big|_{a=1-}-\partial_a^{+} \ent_{\vol_{g_\lambda}}(\mu^a_\lambda)\big|_{a=0+}\\
		\geq &\int_{U_0}\int_0^1 | S_{\lambda,x}(a)|_{\HS}^2\d a\d\mu^0_\lambda(x)+K\lambda^2W^2_{2,g_{\lambda}}(\mu_\lambda^0,\mu_\lambda^1)\geq c_0+K\lambda^2W^2_{2,g_{\lambda}}(\mu_\lambda^0,\mu_\lambda^1)\,,\label{ineq:uniform-eta-lower}
	\end{align}
	where we used \eqref{ineq:Ric_g_lambda-2}, \eqref{ineq:uniform-Riccati-coercivity}, and the fact that $ S_{\lambda,x}$ defined by \eqref{eq:tensorS} solves \eqref{eq:Riccati_proof1} with $y=T_\lambda(x)$ for $\mu_\lambda^0$-a.e. $x$, where $T_\lambda$ denotes the optimal map.\medskip
	\\
	\textbf{Step 2}: Scaling back to the original metric.
	By definition, the map
	\begin{equation}
		\iota_\lambda\colon\bigl((0,l/\lambda)\times F,\iota_\lambda^*g\bigr)\to (U_o\setminus\{o\},g)
	\end{equation}
	is an isometry.
	Furthermore, since $g_{\lambda}$ is a constant scaling of $\iota_\lambda^* g$, a curve $\gamma$ in $(0,l/\lambda)\times F$ is a $g_\lambda$-geodesic if and only if $\iota_\lambda\circ\gamma$ is a $g$-geodesic in $U_o\setminus\{o\}$, and
	\begin{equation}\label{eq:length_dilation}
		{\rm L}_g(\iota_\lambda \circ \gamma)={\rm L}_{\iota^*_\lambda g}(\gamma)=\lambda {\rm L}_{g_{\lambda}}(\gamma).
	\end{equation}
	Therefore, by Step~1, $\iota_\lambda(\Omega)$ is a $g$-geodesic tube of $\iota_\lambda(U_0)$ and $\iota_\lambda(U_1)$ for all $\lambda<\lambda_0$.
	
	Now, let $(\mu^a_\lambda)_{a\in[0,1]}$ be a $W_{2,g_\lambda}$-geodesic with endpoints supported in $U_0,U_1$.
	By \eqref{eq:length_dilation}, we have
	\begin{equation}\label{eq:Wasserstein_scaling}
		W_{2,g}((\iota_\lambda)_\#\mu^a_\lambda,(\iota_\lambda)_\#\mu^b_\lambda)
		=
		\lambda W_{2,g_\lambda}(\mu^a_\lambda,\mu^b_\lambda)\qquad \forall a,b\in [0,1].
	\end{equation}
	Thus $\left((\iota_\lambda)_\#\mu^a_\lambda\right)_a$ is a $W_{2,g}$-geodesic.
	Conversely, every $W_{2,g}$-geodesic between measures supported in $\iota_\lambda(U_0)$ and $\iota_\lambda(U_1)$ pulls back to a $W_{2,g_\lambda}$-geodesic via $\iota_\lambda$.
	
	Similar to the proof of \cref{thm:rigidity_E_cone}, since $\iota_\lambda^*\vol_g=\lambda^n\vol_{g_\lambda}$ is a constant multiple, the two entropies differ by the constant $n\log\lambda$, and hence
	\begin{equation}\label{eq:entropy_scaling}
		\partial^{\pm}_a \ent_{\vol_{g_\lambda}}(\mu^a_\lambda)= \partial^{\pm}_a\ent_{\vol_g}((\iota_\lambda)_\#\mu^a_\lambda),\qquad a\in[0,1].
	\end{equation}
	For sufficiently small radii in the definition of $\eta(\iota_\lambda(p),\iota_\lambda(q))$, all the corresponding optimal transports remain in the geodesic tube $\iota_\lambda(\Omega)$.
	Therefore, by \eqref{ineq:uniform-eta-lower}, \eqref{eq:Wasserstein_scaling}, and \eqref{eq:entropy_scaling},
	\begin{align}
		\eta_g(\iota_\lambda(p),\iota_\lambda(q))&\geq \inf_{\left((\iota_\lambda)_\#\mu^a_\lambda\right)} \frac{1}{W^2_{2,g}((\iota_\lambda)_\#\mu^0_\lambda,(\iota_\lambda)_\#\mu^1_\lambda)}\mathcal{DE}_{\vol_{g}}((\iota_\lambda)_\#\mu^{\boldsymbol{\cdot}}_{\lambda})\\
		&= \inf_{(\mu^a_\lambda)}\frac{1}{\lambda^2 W^2_{2,g_\lambda}(\mu^0_\lambda,\mu^1_\lambda)}\mathcal{DE}_{\vol_{g_\lambda}}(\mu^{\boldsymbol{\cdot}}_{\lambda})\\
		&\geq K+c_0\, \inf_{(\mu^a_\lambda)}\left(\lambda^2 W^2_{2,g_\lambda}(\mu^0_\lambda,\mu^1_\lambda)\right)^{-1}
	\end{align}
	where $(\mu^a_\lambda)_a$ is taken over all $W_{2,g_\lambda}$-geodesics whose endpoints have finite entropy and are supported in $U_0$ and $U_1$, respectively.
	Since $\iota_\lambda (p),\iota_\lambda (q)\to o$ as $\lambda\to 0$, this implies $\eta^*(o)=\infty$, contradicting the assumption.
	In conclusion, $\Sec_{h_0}\equiv1$.
	Finally, the equivalence on the removability of the conic singularity is given by \cref{lemma:removable_singularity}.
\end{proof}

\begin{corollary}
	\cref{thm:intro-exact-cone} holds for every $K\in\R$.
\end{corollary}
\begin{proof}
	The proof of \cref{thm:rigidity_E_cone} works for general $K\in\R$ except for the derivation of $\Sec_h\equiv 1$ in \ref{item:Econe-rigidity2}$\Rightarrow$\ref{item:Econe-rigidity3}, at which stage $N=n$ and the constancy of $f$ are already established.
	Thus by \cref{thm:intro-conic},\footnote{Although \cref{def:conic-singularity} requires $F$ to be boundary-free, the proof of $\Sec_{h_0}\equiv1$ in \cref{thm:intro-conic} only uses a geodesic tube in $(0,\infty)\times F^\circ$ and the local convergence of $g_\lambda$, and hence applies here. Note that the boundary-free condition is mainly used for \cref{prop:equiv_lower_bounds_conic} and \cref{lemma:removable_singularity}.} it suffices to show that the $(K,n)$-cone $\Sigma$ over $(F,h)$, regarded as an $(n+1)$-dimensional Riemannian manifold with a conic singularity at $o$, is of order $2$, with
	\begin{equation}
		h_r=\frac{\sin_K^2(r)}{r^2}h,\qquad h_0=h.
	\end{equation}
	We have
	\begin{equation}
		g-g_\Sigma=(\sin_K^2(r)-r^2)h
		=a(r)(g_\Sigma-\d r\otimes\d r),
		\qquad
		a(r)\coloneqq\frac{\sin_K^2(r)-r^2}{r^2}.
	\end{equation}
	The Taylor expansion of $\sin_K$ at $0$ gives
	\begin{equation}
		a(r)=-\frac K3r^2+O(r^4),\qquad
		a'(r)=O(r),\qquad
		a''(r)=O(1).
	\end{equation}
	By the Hessian formula for the radial function on a warped product
	\begin{equation}
		\nabla^{g_\Sigma}\d r
		={\rm Hess}_{g_\Sigma}r
		=\frac1r(g_\Sigma-\d r\otimes\d r).
	\end{equation}
	Since $\nabla^{g_\Sigma}g_\Sigma=0$, differentiating this identity and applying the Leibniz rule yield
	\begin{equation}
		\left|(\nabla^{g_\Sigma})^j(g_\Sigma-\d r\otimes\d r)\right|_{g_\Sigma}
		=O(r^{-j}),\qquad 0\leq j\leq2,
	\end{equation}
	where we used $|g_\Sigma-\d r\otimes\d r|_{g_\Sigma}=\sqrt{n}$.
	Moreover,
	\begin{equation}
		(\nabla^{g_\Sigma})^2a
		=a''(r)\d r\otimes\d r+a'(r)\nabla^{g_\Sigma}\d r,
	\end{equation}
	and hence
	\begin{equation}
		\left|(\nabla^{g_\Sigma})^ja\right|_{g_\Sigma}
		=O(r^{2-j}),\qquad 0\leq j\leq2.
	\end{equation}
	Another application of the Leibniz rule therefore gives
	\begin{equation}
		\left|(\nabla^{g_\Sigma})^j(g-g_\Sigma)\right|_{g_\Sigma}
		=O(r^{2-j}),\qquad 0\leq j\leq2.
	\end{equation}
	Consequently, \eqref{eq:asymptotic_cone} holds with $k=2$ and every $\alpha\in(0,2)$.
\end{proof}
\subsection{Applications}\label{sec:application}
In this subsection, we give two applications of the cone rigidity theorem to synthetic upper Ricci bounds and synthetic Ricci flows.
\subsubsection*{Instability of synthetic Ricci upper bounds}
The Stenzel metric \cite{Stenzel} is a smooth, complete Ricci-flat K\"ahler metric $g$ on $T^*\mathbb{S}^n$, the cotangent bundle of the $n$-sphere, for $n\geq 2$; here $T^*\mathbb S^n$ is identified with the affine quadric $\{z\in\mathbb C^{n+1}:\sum_{i=1}^{n+1}z_i^2=1\}$.
For $n=2$, it is the Eguchi--Hanson metric \cite{Eguchi-Hanson79}, whose tangent cone at infinity $\mathbb C^2/\mathbb Z_2$ is flat; see \cite[Example~4.1 and Remark~5.10]{Conlon-Hein13}.
For $n\geq3$, it is asymptotically conical in the sense of \cite[Definition~1.11]{Conlon-Hein13}, and its tangent cone at infinity is the ordinary double point $C\coloneqq\{z\in\mathbb C^{n+1}:\sum_{i=1}^{n+1}z_i^2=0\}$ endowed with its Ricci-flat K\"ahler cone metric; see \cite[Sections~5.4 and~5.6]{Conlon-Hein13}.
 
We consider the case $n=3$.
Then $C\subset\mathbb C^4$ is the conifold, and, as a Riemannian cone, it is the Euclidean cone over the Stiefel manifold $V_2(\mathbb R^4)\cong\mathbb S^2\times\mathbb S^3$ endowed with the homogeneous Sasaki--Einstein metric $T^{1,1}$, which is regular over $\mathbb{CP}^1\times\mathbb{CP}^1$; see \cite[Section~2.3]{Sparks11}.
Let $o$ denote its vertex.
For any base point $p$ and any sequence $\lambda_i\downarrow0$, the blow-down sequence $(T^*\mathbb S^3,\lambda_i^2g,p)$ converges to $(C,o)$ in the pointed Gromov--Hausdorff sense (see \cite[Section~1.1]{Conlon-Hein13}), and, by the asymptotically conical condition, in the smooth Cheeger--Gromov sense on compact subsets of $C\setminus\{o\}$.
Moreover, $C$ is not flat, i.e.\ $\Sec_{T^{1,1}}\not\equiv1$: since $\mathbb{S}^2\times\mathbb{S}^3$ is simply connected but not homeomorphic to $\mathbb S^5$, it admits no metric of constant sectional curvature $1$ by the Killing--Hopf theorem.
 
Since $T^{1,1}$ is Einstein with $\Ric=4g_{T^{1,1}}$, the cone $C$, endowed with the Riemannian volume measure, satisfies $\rcd(0,6)$ by \cite{Ketterer-JMPA15}; in particular, $\eta_*(o)\geq0$.
Hence \cref{thm:rigidity_E_cone}, applied with $F=T^{1,1}$, yields $\eta^*(o)=+\infty$.
On the other hand, every space $(T^*\mathbb S^3,\lambda_i^2g)$ in the blow-down sequence is smooth, complete, Ricci-flat and without boundary, and thus satisfies $\eta_*\equiv\eta^*\equiv 0$.
Thus, the synthetic upper Ricci bound formulated via $\eta^*$, as considered in \cite{Sturm2017remarks,Erbar-Li-Schultz26}, is not stable under Gromov--Hausdorff convergence.
 
It is worth mentioning that \cite[Section~5]{Sturm2017remarks} introduced a stable notion of synthetic Ricci upper bounds by additionally requiring a uniform extendability of geodesics; see also Chapter~5.5 of the author's thesis \cite{Li26} for its generalization to the dynamic setting.
In particular, the blow-down sequence above cannot satisfy this extendability condition uniformly.

\subsubsection*{Comparison between different notions of synthetic Ricci flows}
    In \cite{Erbar-Li-Schultz26}, two notions of Ricci flow for time-dependent metric measure spaces were proposed.
    The first uses dynamic convexity and local almost concavity of the entropy, which can be equivalently formulated via time-dependent $\eta$-quantities.
    The second uses contraction estimates and local short-time expansion estimates for the transport cost between heat flows.
    These are called \textbf{weak Ricci flows} and \textbf{rough Ricci flows}, respectively; see
    \cite[Section~3]{Erbar-Li-Schultz26} for the precise definitions.
    The term ``weak'' is motivated by the heuristic that the entropic formulation is weaker, in the sense that every rough Ricci flow is a weak Ricci flow (see \cite[Corollary~5.8]{Erbar-Li-Schultz26}), whereas the converse fails in general.
    
    In \cite{Erbar-Li-Schultz26}, the failure of equivalence between weak and rough Ricci flows was rigorously demonstrated, but only for two-dimensional cones.
    The argument uses the fact that the regular parts of Euclidean and spherical cones are locally isometric to open subsets of $\R^2$ and $\S^2$, respectively.
    Our cone rigidity result \cref{thm:intro-exact-cone}, combined with \cite[Theorem~6.1]{Erbar-Li-Schultz26}, extends this comparison to arbitrary dimensions and yields the following characterization.
    \begin{theorem}
     Let $(F,h)$ be a closed $n$-dimensional smooth Riemannian manifold, with $n\geq 2$ and $\Ric_h=(n-1)h$.
      Consider the $(1,n)$-cone $\Sigma_1(F)$, with cone distance $d_{\Sigma_1}$ and Riemannian volume measure $\m_{\Sigma_1^n}$, as in \cref{sec:cones}.
      Define
    	\begin{equation}
    		d_t\coloneqq (1-2nt)^{\frac12}d_{\Sigma_1},\qquad \m_t\coloneqq (1-2nt)^{\frac{n+1}{2}}\m_{\Sigma^n_1}.
    	\end{equation}
    Then the time-dependent metric measure space $\left( \Sigma_1(F),d_t,\m_t\right)_{t\in[0,\frac{1}{2n})}$ is
    \begin{itemize}
    	\item a weak Ricci flow in the sense of \cite[Definition 3.1]{Erbar-Li-Schultz26} if and only if $\Sec_h\equiv 1$ on $F$;
    	\item a rough Ricci flow in the sense of \cite[Definition 3.7]{Erbar-Li-Schultz26} if and only if $(F,h)$ is isometric to the unit sphere $\mathbb{S}^n$.
    \end{itemize}
    \end{theorem}
In particular, for $F=\mathbb{RP}^4=\S^4/\mathbb{Z}_2$, the time-scaled spherical cone $(\Sigma_1(F),d_t,\m_t)_{t\in[0,\frac{1}{8})}$ is a weak Ricci flow but not a rough Ricci flow; and for $F=\mathbb{S}^2(1/\sqrt{3})\times \mathbb{S}^2(1/\sqrt{3})$, the corresponding time-scaled spherical cone is neither a weak Ricci flow nor a rough Ricci flow.
    \appendix
\section{Proof of the boundary transport lemma}\label{app:boundary-transport}
\begin{proof}[Proof of \cref{lem:boundary-transport}]
	Choose a smooth hypersurface $S\subset\partial M$ through $p$ whose unit normal at $p$ is $v$, and let $\rho$ be its signed distance in the boundary metric, with sign chosen so that $\nabla^{\partial M}\rho(p)=v$.
	Denote $v_q\coloneqq\nabla^{\partial M}\rho(q)$.
	After shrinking $\Sigma$, we may assume that
	\begin{equation}\label{eq:II<0}
		\lambda_q\coloneqq-\II_q(v_q,v_q)\geq\lambda_0>0,
		\qquad \text{on $\Sigma$}.
	\end{equation}
	
	\textbf{Part 1.} \cref{item:bdy_calculus1}.\smallskip
	
	Let $(\widetilde M,\widetilde g)$ be a smooth extension of $M$ across $\Sigma$.
	Consider the signed geodesic map
	\begin{equation}
		E\colon\Sigma\times(-r_0,r_0)\to\widetilde M,
		\qquad
		E(q,s)\coloneqq\exp_q^{\widetilde M}(s v_q),
	\end{equation}
	so that $E_\pm(q,t)=E(q,\pm t)$ for $t\in(0,r_0)$.
	Let $r$ be the normal signed distance to the boundary.
	Since $\partial_sE(q,0)=v_q$, \eqref{eq:II-and-Hess} gives
	\begin{align}
		\partial_s\bigl(r(E(q,s))\bigr)|_{s=0}
		&=\widetilde g_q(\nabla^{\widetilde g}r,v_q)=0,\\
		\partial_s^2\bigl(r(E(q,s))\bigr)|_{s=0}
		&=\Hess_{\widetilde g}r(v_q,v_q)=-\II_q(v_q,v_q)=\lambda_q.
	\end{align}
	Therefore,
	\begin{equation}\label{eq:r-Phi-expansion}
		r(E(q,s))
		=\frac12\lambda_qs^2+O(s^3),
		\qquad
		\partial_s\bigl(r(E(q,s))\bigr)
		=\lambda_qs+O(s^2),
	\end{equation}
	uniformly in $q$.
	In particular, $E_\pm(q,t)\in M^\circ$ for all sufficiently small $t>0$.
	
	Since $E=\mathrm{id}$ on $\Sigma\times\{0\}$, the differential of $E$ on the boundary satisfies
	\begin{equation}\label{eq:dPhi-boundary}
		\d E_{(q,0)}(w,\tau)=w+\tau v_q,
		\qquad w\in T_q\partial M,\ \tau\in\R.
	\end{equation}
	Although $\d E$ is degenerate on $\Sigma\times\{0\}$ as $\d E_{(q,0)}(-v_q,1)=0$, the maps $E_\pm$ are diffeomorphisms on $\Sigma\times(0,r_0)$ onto their images.
	To see this, Taylor's formula gives a smooth positive function $b$, after shrinking the neighborhoods, such that
	\begin{equation}
		r(E(q,s))=\frac12s^2b(q,s),
		\qquad b(q,0)=\lambda_q.
	\end{equation}
	Define $R(q,s)\coloneqq s\sqrt{b(q,s)}$ for $-r_0<s<r_0$.
	Then $R$ is smooth across $s=0$ and $r(E(q,s))=R(q,s)^2/2$.
	
	Let $P$ denote the collar projection, i.e.\ $P({\mathcal C}(z,r))=z$, where ${\mathcal C}$ is the collar map in \eqref{eq:normal_collar_map}.
	Define $\Psi\colon\Sigma\times(-r_0,r_0)\to\partial M\times\R$ and $H$ by
	\begin{equation}
		\Psi(q,s)\coloneqq\bigl(P(E(q,s)),R(q,s)\bigr),
		\qquad H(z,s)\coloneqq(z,s^2/2).
	\end{equation}
	By construction,
	\begin{equation}\label{eq:Phi-composition}
		E={\mathcal C}\circ H\circ\Psi.
	\end{equation}
	Since $\d P_q$ is the identity on $T_q\partial M$ and $\partial_sR(q,0)=\sqrt{\lambda_q}$, we have
	\begin{equation}\label{eq:dPsi}
		\d\Psi_{(q,0)}(w,\tau)
		=\bigl(w+\tau v_q,\tau\sqrt{\lambda_q}\bigr).
	\end{equation}
	Applying the inverse function theorem at $(p,0)$, we may shrink $\Sigma$ and $r_0$ so that $\Psi$ is a diffeomorphism on $\Sigma\times(-r_0,r_0)$ onto its image.
	Since $R(q,s)$ has the same sign as $s$, and $H$ is a diffeomorphism on each side of $\{s=0\}$, \eqref{eq:Phi-composition} gives the same conclusion for $E_+$ and $E_-$.
	Both images contain the interior of a common collar neighborhood $\Omega$ of $p$ in $M$.
	For sufficiently small $\delta$, $\Omega$ contains the boundary geodesic from $x$ to $y$.
	This proves \cref{item:bdy_calculus1}.
	\medskip
	
	\textbf{Part 2.} The gradient identities in \cref{item:bdy_calculus2}.\smallskip
	
	Define $u_\pm$ on $\Omega\cap M^\circ$ by \eqref{eq:u-pm}, or equivalently,
	\begin{equation}\label{eq:u_pm_appendix}
		u_\pm\bigl(E(q,s)\bigr)=\rho(q)+s,
		\qquad \pm s\in(0,r_0).
	\end{equation}
	By \cref{item:bdy_calculus1}, $u_\pm$ are smooth.
	Since the two branches of $E$ (i.e.\ the restrictions of $E$ to $\{\pm s>0\}$, respectively) extend continuously to $s=0$, both functions extend continuously to $\Omega\cap\partial M$, with $u_\pm|_{\Omega\cap\partial M}=\rho$.
	
	On either branch, differentiating \eqref{eq:u_pm_appendix} yields
	\begin{equation}\label{eq:variation-field-Phi}
		\d u_\pm\bigl(\d E_{(q,s)}(w,0)\bigr)=\d\rho_q(w),
		\qquad
		\d u_\pm\bigl(\partial_sE(q,s)\bigr)=1,
	\end{equation}
	where the sign of $s$ determines the corresponding branch.
	For $w\in T_q\partial M$, the standard Jacobi field computation for the unit-speed geodesic variation $E$ gives
	\begin{equation}\label{eq:Gauss-signed-Phi}
		\widetilde g\bigl(\partial_sE(q,s),\d E_{(q,s)}(w,0)\bigr)
		=\widetilde g(v_q,w)=\d\rho_q(w).
	\end{equation}
	Since $\d E_{(q,s)}$ is invertible for $s\neq0$, comparing \eqref{eq:variation-field-Phi} and \eqref{eq:Gauss-signed-Phi} gives
	\begin{equation}\label{eq:gradient-u-signed}
		\nabla u_\pm(E(q,s))=\partial_sE(q,s),
		\qquad \pm s\in(0,r_0).
	\end{equation}
	As $\partial_tE_\pm(q,t)=\pm\partial_sE(q,\pm t)$, this shows \eqref{eq:integral-curve-u_pm}, and by the definition of $E$ we get $|\nabla u_\pm|=1$.\medskip
	
	\textbf{Part 3.} The Laplacian estimates for $u_\pm$.\smallskip
	
	It suffices to consider $\Delta u_+$; the negative branch is treated in the same way.
	Let $\sigma$ denote the Riemannian volume measure on $\partial M$ and $J$ the Jacobian of $E$ on either branch s.t.
	\begin{equation}
		E^*\vol_{\widetilde g}=J(q,s)\d\sigma(q)\d s,
		\qquad J(q,s)=\det\d E_{(q,s)}.
	\end{equation}
	
	Let $\chi\in C_c^\infty(\Omega\cap M^\circ)$.
	By \eqref{eq:gradient-u-signed}, the chain rule, and change of variables, we have
	\begin{align}
		\int_{\Omega\cap M^\circ}\chi\Delta u_+\d\vol_{\widetilde g}
		&=-\int_{\Omega\cap M^\circ}
		\widetilde g(\nabla\chi,\nabla u_+)\d\vol_{\widetilde g}\\
		&=-\int\partial_s(\chi\circ E)(q,s)
		J(q,s)\d\sigma(q)\d s\\
		&=\int(\chi\circ E)(q,s)
		\partial_sJ(q,s)\d\sigma(q)\d s.
	\end{align}
	As $\chi$ is arbitrary,
	\begin{equation}\label{eq:laplacian-u-J}
		\Delta u_+(E(q,s))
		=\frac{\partial_sJ(q,s)}{J(q,s)}.
	\end{equation}
	Among the components of $E$ in \eqref{eq:Phi-composition}, the Jacobian of $H$ at $\Psi(q,s)$ is $|R(q,s)|$, and the Jacobians for ${\mathcal C}$ and $\Psi$ are smooth and positive near $s=0$.
	Thus $J(q,s)=|s|\,\beta(q,s)$ for a smooth positive function $\beta$ extending across $s=0$.
	After shrinking the neighborhoods, $\partial_s\beta/\beta$ is uniformly bounded.
	Consequently, on the respective branches,
	\begin{equation}
		\Delta u_\pm(E(q,s))
		=\partial_s\log J(q,s)=\frac1s+O(1).
	\end{equation}
	Substituting $s=\pm t$ proves \eqref{eq:laplacian-estimates}, uniformly in $q$.\medskip
	
	\textbf{Part 4.} \cref{item:bdy_calculus3}.\smallskip
	
	Let $x,y$ be as in \eqref{eq:xy}, with $\delta>0$ sufficiently small, and let $\bar\gamma\colon[0,1]\to\widetilde M$ be the $\widetilde g$-geodesic from $x$ to $y$.
	As $\delta\downarrow0$, its image converges to $p$, while the direction of its velocity converges uniformly to $v$.
	Thus by \eqref{eq:II<0}, for sufficiently small $\delta$, we have
	\[
	(r\circ\bar\gamma)''=\Hess_{\widetilde g}r\,(\bar\gamma',\bar\gamma')>0.
	\]
	As $x,y\in\partial M$, $r(\bar\gamma(0))=r(\bar\gamma(1))=0$.
	Strict convexity gives $r(\bar\gamma(a))<0$ for all $a\in(0,1)$.
	By smooth dependence on the endpoints, there are neighborhoods $U_x,U_y\subset\Omega$ of $x,y$ s.t. the minimizing $\widetilde g$-geodesic from any $\tilde x\in U_x\cap M^\circ$ to $\tilde y\in U_y\cap M^\circ$ also meets $\{r<0\}$, and hence intersects $\partial M$.
	
	Since $|\nabla u_\pm|=1$, every rectifiable curve $\zeta$ in $\Omega\cap M^\circ$ satisfies
	\begin{equation}
		u_\pm(\zeta(1))-u_\pm(\zeta(0))\leq\mathrm L_g(\zeta),
	\end{equation}
	where $\mathrm L_g$ stands for the length w.r.t.\ $g$.
	The same inequality holds for curves meeting $\partial M$, by applying it to inward collar shifts whose lengths converge to the original length, and using the continuity of $u_\pm$.
	
	Now, assume $\gamma\colon[0,1]\to M$ is a minimizing geodesic from $\tilde x$ to $\tilde y$ w.r.t. $d_g$.
	From the previous discussion, one can see that $\gamma$ has to intersect $\partial M$ at a point $\gamma(a)$ with $a\in(0,1)$.
	Then, as $u_-(\gamma(a))=u_+(\gamma(a))=\rho(\gamma(a))$, it follows that
	\begin{align}
		d_g(\tilde x,\tilde y)
		&=\mathrm L_g(\gamma|_{[0,a]})+\mathrm L_g(\gamma|_{[a,1]})\\
		&\geq u_-(\gamma(a))-u_-(\tilde x)+u_+(\tilde y)-u_+(\gamma(a))=u_+(\tilde y)-u_-(\tilde x).
	\end{align}
	This proves \cref{item:bdy_calculus3}.
\end{proof}

   \section{Semi-convexity of the squared distance between geodesics}
  We present a semi-convexity estimating for the squared distance between two geodesics (not necessarily emanate from a common point) on ${\rm CAT}(\kappa)$ spaces.
  Such a result may be of independent interest.
  For $\kappa=0$, it reduces to the classical convexity of the metric in $\mathrm{CAT}(0)$ spaces; see \cite[Proposition~II.2.2]{Bridson-Haefliger99} or \cite[Corollary 2.5]{Sturm03}.
 
  Set $\beta_0\equiv0$ and, for $\kappa>0$, define
  \begin{equation}
  	\beta_\kappa(r)\coloneqq
  	\sqrt\kappa\,r\tan\bigl(\sqrt\kappa\,r/2\bigr),
  	\qquad 0\leq r<\pi/\sqrt\kappa.
  \end{equation}
  \begin{lemma}\label{lemma:convexityCAT}
  	Let $\kappa\geq0$ and let $(U,d)$ be a $\mathrm{CAT}(\kappa)$ space.
  	Let $\gamma,\widetilde\gamma\colon I\to U$ be constant-speed
  	geodesics on an open interval $I$.
  	If $\kappa>0$, assume that
  	$d(\gamma^a,\widetilde\gamma^a)<\pi/\sqrt\kappa$ for all $a\in I$.
  	Then
  	\begin{equation}\label{eq:CAT-joint-comparison}
  		\partial_a^2\left(\frac12d^2(\gamma^a,\widetilde\gamma^a)\right)
  		\geq
  		-\beta_\kappa\bigl(d(\gamma^a,\widetilde\gamma^a)\bigr)
  		\bigl(|\dot\gamma|^2+|\dot{\widetilde\gamma}|^2\bigr)
  	\end{equation}
  	holds in the distributional sense on $I$.
  \end{lemma}
  
  \begin{proof}
  	For $\kappa=0$ the claim follows from the convexity of
  	$a\mapsto d(\gamma^a,\widetilde\gamma^a)$,
  	see \cite[Proposition~II.2.2]{Bridson-Haefliger99}.
  	Assume $\kappa>0$.
  	Let $(\Sigma,d_\Sigma)$ be the Euclidean cone over
  	$(U,\sqrt\kappa\,d)$, which is $\mathrm{CAT}(0)$ by
  	\cite[Theorem~II.3.14]{Bridson-Haefliger99}.
  	Consider the lifts of $\gamma$ and $\widetilde\gamma$ to
  	$\{1\}\times U$; by the definition of the cone distance in \cref{def:cone}, 
  	\begin{equation}\label{eq:appendix_u}
  		L(a)\coloneqq\frac12d_\Sigma^2
  		\bigl((1,\gamma^a),(1,\widetilde\gamma^a)\bigr)
  		=
  		1-\cos\bigl(\sqrt\kappa\,d(\gamma^a,\widetilde\gamma^a)\bigr).
  	\end{equation}
  	Denote $\ell=|\dot\gamma|$ and $\widetilde\ell=|\dot{\widetilde\gamma}|$.
  	Fix $a$ and take $h>0$ sufficiently small.
  	In the rescaled base $(U,\sqrt\kappa\,d)$, the restriction
  	$\gamma|_{[a-h,a+h]}$ is a geodesic segment of length
  	$2\sqrt\kappa\,\ell h<\pi$.
  	Note that the cone over the geodesic segment is a Euclidean sector.
  	Thus the midpoint of $(1,\gamma^{a-h})$ and $(1,\gamma^{a+h})$ with
  	respect to the cone metric is $(r,\gamma^a)$, where
  	$r=\cos(\sqrt\kappa\,\ell h)$; see \cref{fig:cone-midpoint}.
  	Similarly, the cone midpoint of $(1,\widetilde\gamma^{a-h})$ and
  	$(1,\widetilde\gamma^{a+h})$ is $(\widetilde r,\widetilde\gamma^a)$ with
  	$\widetilde r=\cos(\sqrt\kappa\,\widetilde\ell h)$.
  	\begin{figure}[htbp]
  		\centering
        \begin{tikzpicture}[scale=5.2, line cap=round, line join=round,
    font=\small]
    \definecolor{coneaccent}{RGB}{35,78,113}
    \def\SectorAngle{34}
    \coordinate (O) at (0,0);
    \coordinate (L) at ({-sin(\SectorAngle)},{-cos(\SectorAngle)});
    \coordinate (R) at ({ sin(\SectorAngle)},{-cos(\SectorAngle)});
    \coordinate (M) at (0,{-cos(\SectorAngle)});
    \coordinate (Q) at (0,-1);

    \fill[black!2] (O) -- (L)
        arc[start angle={-90-\SectorAngle},
            end angle={-90+\SectorAngle}, radius=1] -- cycle;
    \draw[line width=0.7pt] (L) -- (O) -- (R);
    \draw[black!55, densely dashed, line width=0.6pt] (L)
        arc[start angle={-90-\SectorAngle},
            end angle={-90+\SectorAngle}, radius=1];

    \draw[coneaccent, line width=1pt] (L) -- (R);
    \draw[black!70, line width=0.65pt] (O) -- (M);
    \draw[black!55, densely dashed, line width=0.6pt] (M) -- (Q);

    \draw[black!65, line width=0.5pt]
        ({-90-\SectorAngle}:0.22)
        arc[start angle={-90-\SectorAngle}, end angle=-90, radius=0.22];
    \draw[black!65, line width=0.5pt] (-90:0.22)
        arc[start angle=-90, end angle={-90+\SectorAngle}, radius=0.22];
    \node at ({-90-\SectorAngle/2}:0.32) {$\theta$};
    \node at ({-90+\SectorAngle/2}:0.32) {$\theta$};
    \draw[black!65, line width=0.45pt]
        (0,{-cos(\SectorAngle)+0.055}) --
        (0.055,{-cos(\SectorAngle)+0.055}) --
        (0.055,{-cos(\SectorAngle)});

    \path (O) -- (L) node[pos=0.54, above left=2pt] {$1$};
    \path (O) -- (R) node[pos=0.54, above right=2pt] {$1$};
    \path (O) -- (M) node[pos=0.67, left=5pt] {$r$};

    \fill (O) circle[radius=0.3pt] node[above=5pt] {$o$};
    \fill (L) circle[radius=0.3pt]
        node[left=5pt] {$(1,\gamma^{a-h})$};
    \fill (R) circle[radius=0.3pt]
        node[right=5pt] {$(1,\gamma^{a+h})$};
    \fill[coneaccent] (M) circle[radius=0.34pt]
        node[above right=8pt] {$(r,\gamma^a)$};
    \draw[fill=white, black!70, line width=0.6pt]
        (Q) circle[radius=0.3pt] node[below=6pt] {$(1,\gamma^a)$};
     \node at (0,-1.28)
        {$\theta=\sqrt\kappa\,\ell h,\qquad r=\cos\theta$};
\end{tikzpicture}
  		\caption{The subcone over $\gamma([a-h,a+h])$ is a Euclidean
  			sector. The chord joining the two unit-radius endpoints has
  			midpoint $(r,\gamma^a)$, where
  			$r=\cos(\sqrt\kappa\,\ell h)$.
  			The dashed arc lies at unit radius.}
  		\label{fig:cone-midpoint}
  	\end{figure}
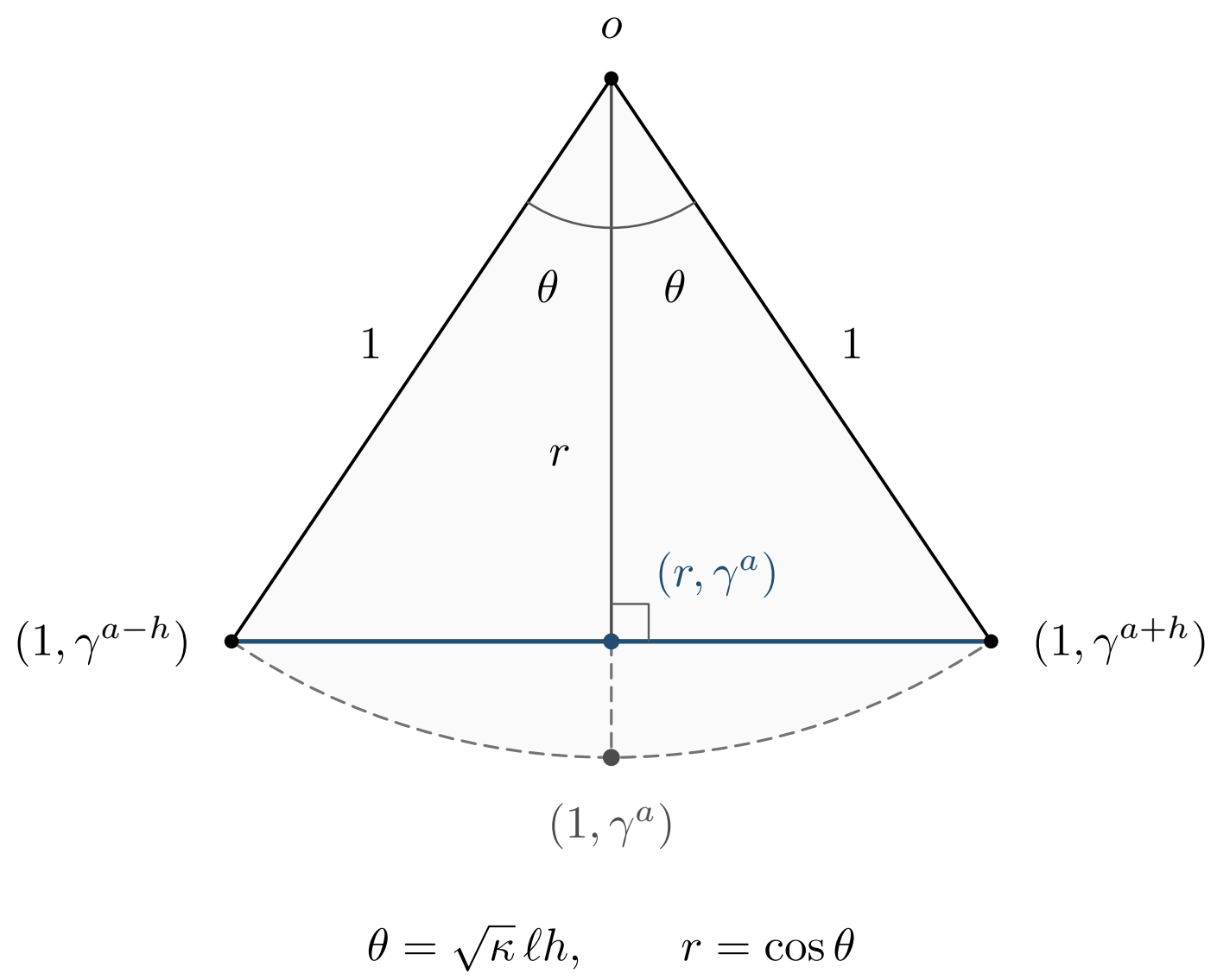
  	
  	\noindent Now applying convexity of the squared distance in the cone,
  	together with \eqref{eq:appendix_u}, yields
  	\begin{align}
  		L(a-h)+L(a+h)&\geq d^2_\Sigma\bigl((r,\gamma^a),(\widetilde r, \widetilde \gamma^a)\bigr)
  		=r^2+\widetilde r^2-2r\widetilde r \cos\bigl(\sqrt\kappa\,d(\gamma^a,\widetilde \gamma^a)\bigr)\\
  		&=(r-\widetilde r)^2+2r\widetilde r L(a).
  	\end{align}
  	Since $0<r,\widetilde r\leq1$ and $1-\cos z\leq z^2/2$, we have
  	$1-r\widetilde r\leq(1-r)+(1-\widetilde r)\leq
  	\kappa(\ell^2+\widetilde\ell^2)h^2/2$.
  	Hence
  	\begin{align}
  	   L(a-h)+L(a+h)-2L(a)&\geq (r-\widetilde r)^2+(2r\widetilde r -2) L(a)\\
       &\geq (2r\widetilde r -2) L(a)\geq-\kappa(\ell^2+\widetilde\ell^2)L(a)h^2.\label{eq:CAT-midpoint-estimate}
  	\end{align}
  	Thus $L''\geq-\kappa(\ell^2+\widetilde\ell^2)L$ in the sense of
  	distributions; in particular $L$ is locally semiconvex, hence locally
  	Lipschitz, and $L''$ is a locally finite measure.
  	Denote by $\phi$ the function on $[0,2)$ given by
  	\begin{equation}
  		\phi(\tau)\coloneqq \frac{[\arccos(1-\tau)]^2}{2\kappa},\qquad
  		\phi'(\tau)=\frac{\arccos(1-\tau)}{\kappa\sqrt{1-(1-\tau)^2}}.
  	\end{equation}
  	Then direct computation gives
  	\begin{equation}\label{eq:appendix_Phi}
  		\phi(L(a))=\frac12d^2(\gamma^a,\widetilde\gamma^a),\qquad
  		\kappa L(a)\phi'(L(a))
  		=\beta_\kappa\bigl(d(\gamma^a,\widetilde\gamma^a)\bigr).
  	\end{equation}
  	Since $\phi\in C^\infty([0,2))$ with $\phi'\geq0$ and $\phi''\geq0$, the
  	chain rule (valid for the locally semiconvex function $L$, with all terms
  	understood as measures) together with \eqref{eq:CAT-midpoint-estimate}
  	and \eqref{eq:appendix_Phi} gives
  	\begin{align}
  		\partial_a^2\left(\frac12d^2(\gamma^a,\widetilde\gamma^a)\right)
  		&=(\phi\circ L)''=\phi''(L)(L')^2+\phi'(L)L''\geq \phi'(L)L''\\
  		&\geq-\kappa(\ell^2+\widetilde\ell^2)L\phi'(L)
  		=-(\ell^2+\widetilde\ell^2)\beta_\kappa\bigl(d(\gamma^a,\widetilde\gamma^a)\bigr).\tag*{\qedhere}
  	\end{align}
  \end{proof}
    
\bibliographystyle{alpha}
\bibliography{ms}

@article{ABB1993,
	author = {Alexander, Stephanie B. and Berg, I. David and Bishop, Richard L.},
	title = {Geometric curvature bounds in {Riemannian} manifolds with boundary},
	fjournal = {Transactions of the American Mathematical Society},
	journal = {Trans. Am. Math. Soc.},
	issn = {0002-9947},
	volume = {339},
	number = {2},
	pages = {703--716},
	year = {1993},
	language = {English},
	doi = {10.2307/2154294},
	zbMATH = {445013},
	Zbl = {0793.53039}
}

@misc{ARS26,
	author = {Gautam Aishwarya and Liran Rotem and Yair Shenfeld},
	title = {Sectional curvature and matrix displacement convexity},
	year = {2026},
	howpublished = {Preprint, {arXiv}:2509.23399 [math.{DG}] (2026)},
	url = {https://arxiv.org/abs/2509.23399},
	arXiv = {arXiv:2509.23399}
}

@article {AGSRCD2014,
    AUTHOR = {Ambrosio, Luigi and Gigli, Nicola and Savar\'{e}, Giuseppe},
     TITLE = {Metric measure spaces with {R}iemannian {R}icci curvature
              bounded from below},
   JOURNAL = {Duke Math. J.},
  FJOURNAL = {Duke Mathematical Journal},
    VOLUME = {163},
      YEAR = {2014},
    NUMBER = {7},
     PAGES = {1405--1490},
      ISSN = {0012-7094},
   MRCLASS = {35R01 (60J45 60J65)},
  MRNUMBER = {3205729},
       DOI = {10.1215/00127094-2681605},
       URL = {https://doi.org/10.1215/00127094-2681605},
}

@Article{AGS15,
	Author = {Ambrosio, Luigi and Gigli, Nicola and Savar\'{e}, Giuseppe},
	Title = {Bakry--\'{E}mery curvature-dimension condition and {Riemannian} {Ricci} curvature bounds},
	FJournal = {The Annals of Probability},
	Journal = {Ann. Probab.},
	ISSN = {0091-1798},
	Volume = {43},
	Number = {1},
	Pages = {339--404},
	Year = {2015},
	DOI = {10.1214/14-AOP907},
	zbMATH = {6383628},
	Zbl = {1307.49044}
}

@book{Bridson-Haefliger99,
	author = {Bridson, Martin R. and Haefliger, Andr{\'e}},
	title = {Metric spaces of non-positive curvature},
	fseries = {Grundlehren der Mathematischen Wissenschaften},
	series = {Grundlehren Math. Wiss.},
	issn = {0072-7830},
	volume = {319},
	isbn = {3-540-64324-9},
	year = {1999},
	publisher = {Berlin: Springer},
	language = {English},
	zbMATH = {1385418},
	Zbl = {0988.53001}
}

@InCollection{BacherSturm14,
 Author = {Bacher, Kathrin and Sturm, Karl-Theodor},
 Title = {Ricci bounds for {Euclidean} and spherical cones},
 BookTitle = {Singular phenomena and scaling in mathematical models},
 ISBN = {978-3-319-00785-4; 978-3-319-00786-1},
 Pages = {3--23},
 Year = {2014},
 Publisher = {Cham: Springer},
 Language = {English},
 DOI = {10.1007/978-3-319-00786-1_1},
 zbMATH = {6492091},
 Zbl = {1328.53050}
}

@article{Bertrand-Ketterer-Mondello-Richard21,
	author = {Bertrand, J{\'e}r{\^o}me and Ketterer, Christian and Mondello, Ilaria and Richard, Thomas},
	title = {Stratified spaces and synthetic {Ricci} curvature bounds},
	fjournal = {Annales de l'Institut Fourier},
	journal = {Ann. Inst. Fourier},
	issn = {0373-0956},
	volume = {71},
	number = {1},
	pages = {123--173},
	year = {2021},
	language = {English},
	doi = {10.5802/aif.3393},
	zbMATH = {7362100},
	Zbl = {1481.53051}
}

@article{Bishop-ONeill69,
 author = {Bishop, R. L. and O'Neill, B.},
 title = {Manifolds of negative curvature},
 fjournal = {Transactions of the American Mathematical Society},
 journal = {Trans. Am. Math. Soc.},
 issn = {0002-9947},
 volume = {145},
 pages = {1--49},
 year = {1969},
 language = {English},
 doi = {10.2307/1995057},
 zbMATH = {3304660},
 Zbl = {0191.52002}
}

@article{Conlon-Hein13,
	author = {Conlon, Ronan J. and Hein, Hans-Joachim},
	title = {Asymptotically conical {Calabi}-{Yau} manifolds. {I}},
	fjournal = {Duke Mathematical Journal},
	journal = {Duke Math. J.},
	issn = {0012-7094},
	volume = {162},
	number = {15},
	pages = {2855--2902},
	year = {2013},
	language = {English},
	doi = {10.1215/00127094-2382452},
	zbMATH = {6250912},
	Zbl = {1283.53045}
}

@article{Cordero-Erausquin2001,
    AUTHOR = {Cordero-Erausquin, Dario and McCann, Robert J. and
              Schmuckenschl\"{a}ger, Michael},
     TITLE = {A {R}iemannian interpolation inequality \`a la {B}orell,
              {B}rascamp and {L}ieb},
   JOURNAL = {Invent. Math.},
  FJOURNAL = {Inventiones Mathematicae},
    VOLUME = {146},
      YEAR = {2001},
    NUMBER = {2},
     PAGES = {219--257},
      ISSN = {0020-9910},
   MRCLASS = {58E35 (28C99 60E15)},
  MRNUMBER = {1865396},
MRREVIEWER = {C\'{e}dric Villani},
       DOI = {10.1007/s002220100160},
       URL = {https://doi.org/10.1007/s002220100160},
}

@Article{EKS15,
 Author = {Erbar, Matthias and Kuwada, Kazumasa and Sturm, Karl-Theodor},
 Title = {On the equivalence of the entropic curvature-dimension condition and {Bochner}'s inequality on metric measure spaces},
 FJournal = {Inventiones Mathematicae},
 Journal = {Invent. Math.},
 ISSN = {0020-9910},
 Volume = {201},
 Number = {3},
 Pages = {993--1071},
 Year = {2015},
 Language = {English},
 DOI = {10.1007/s00222-014-0563-7},
 zbMATH = {6483768},
 Zbl = {1329.53059}
}

@article{Erbar-Li-Schultz26,
	author = {Erbar, Matthias and Li, Zhenhao and Schultz, Timo},
	title = {Synthetic notions of {Ricci} flow for metric measure spaces},
	fjournal = {Journal of Functional Analysis},
	journal = {J. Funct. Anal.},
	issn = {0022-1236},
	volume = {291},
	number = {12},
	pages = {73},
	note = {Id/No 111661},
	year = {2026},
	language = {English},
	doi = {10.1016/j.jfa.2026.111661},
	zbMATH = {8268291}
}

@Article{Erbar-Sturm,
 Author = {Erbar, Matthias and Sturm, Karl-Theodor},
 Title = {Rigidity of cones with bounded {Ricci} curvature},
 FJournal = {Journal of the European Mathematical Society (JEMS)},
 Journal = {J. Eur. Math. Soc. (JEMS)},
 ISSN = {1435-9855},
 Volume = {23},
 Number = {1},
 Pages = {219--235},
 Year = {2021},
 Language = {English},
 DOI = {10.4171/JEMS/1010},
 zbMATH = {7328111},
 Zbl = {1478.53078}
}

@article{Eguchi-Hanson79,
 author = {Eguchi, Tohru and Hanson, Andrew J.},
 title = {Self-dual solutions to {Euclidean} gravity},
 fjournal = {Annals of Physics},
 journal = {Ann. Phys.},
 issn = {0003-4916},
 volume = {120},
 pages = {82--106},
 year = {1979},
 language = {English},
 doi = {10.1016/0003-4916(79)90282-3},
 url = {www.osti.gov/biblio/1447072},
 zbMATH = {3636825},
 Zbl = {0409.53020}
}

@Book{GigliRCD,
 Author = {Gigli, Nicola},
 Title = {On the differential structure of metric measure spaces and applications},
 FSeries = {Memoirs of the American Mathematical Society},
 Series = {Mem. Am. Math. Soc.},
 ISSN = {0065-9266},
 Volume = {1113},
 ISBN = {978-1-4704-1420-7; 978-1-4704-2279-0},
 Year = {2015},
 Publisher = {Providence, RI: American Mathematical Society (AMS)},
 Language = {English},
 DOI = {10.1090/memo/1113},
 zbMATH = {6479622},
 Zbl = {1325.53054}
}

@article{Han-20,
	author = {Han, Bang-Xian},
	title = {Measure rigidity of synthetic lower {Ricci} curvature bound on {Riemannian} manifolds},
	fjournal = {Advances in Mathematics},
	journal = {Adv. Math.},
	issn = {0001-8708},
	volume = {373},
	pages = {31},
	note = {Id/No 107327},
	year = {2020},
	language = {English},
	doi = {10.1016/j.aim.2020.107327},
	zbMATH = {7243326},
	Zbl = {1447.30007}
}

@misc{Han-17,
	author = {Han, Bang-Xian},
	title = {Ricci tensor on smooth metric measure space with boundary},
	year = {2017},
	howpublished = {Preprint, {arXiv}:1709.10143 [math.{MG}] (2017)},
	url = {https://arxiv.org/abs/1709.10143},
	arXiv = {arXiv:1709.10143}
}

@Article{Ketterer-JFA13,
 Author = {Ketterer, Christian},
 Title = {Ricci curvature bounds for warped products},
 FJournal = {Journal of Functional Analysis},
 Journal = {J. Funct. Anal.},
 ISSN = {0022-1236},
 Volume = {265},
 Number = {2},
 Pages = {266--299},
 Year = {2013},
 Language = {English},
 DOI = {10.1016/j.jfa.2013.05.008},
 zbMATH = {6259496},
 Zbl = {1284.53072}
}

@Article{Ketterer-JMPA15,
 Author = {Ketterer, Christian},
 Title = {Cones over metric measure spaces and the maximal diameter theorem},
 FJournal = {Journal de Math{\'e}matiques Pures et Appliqu{\'e}es. Neuvi{\`e}me S{\'e}rie},
 Journal = {J. Math. Pures Appl. (9)},
 ISSN = {0021-7824},
 Volume = {103},
 Number = {5},
 Pages = {1228--1275},
 Year = {2015},
 Language = {English},
 DOI = {10.1016/j.matpur.2014.10.011},
 zbMATH = {6438813},
 Zbl = {1317.53064}
}

@misc{Ketterer26-RCD-cone,
	author = {Christian Ketterer},
	title = {Warped products over one-dimensional base spaces and the {RCD} condition},
	year = {2026},
	howpublished = {Preprint, {arXiv}:2506.10809 [math.{DG}] (2026)},
	url = {https://arxiv.org/abs/2506.10809},
	arXiv = {arXiv:2506.10809}
}

@article{Ketterer26,
 author = {Ketterer, Christian},
 title = {Warped products and synthetic lower curvature bounds: an overview},
 fjournal = {Indagationes Mathematicae. New Series},
 journal = {Indag. Math., New Ser.},
 issn = {0019-3577},
 volume = {37},
 number = {3},
 pages = {1014--1041},
 year = {2026},
 language = {English},
 doi = {10.1016/j.indag.2025.09.006},
 zbMATH = {8196494}
}

@book{Lee_Riem18,
 author = {Lee, John M.},
 title = {Introduction to {Riemannian} manifolds},
 edition = {2nd},
 fseries = {Graduate Texts in Mathematics},
 series = {Grad. Texts Math.},
 issn = {0072-5285},
 volume = {176},
 isbn = {978-3-319-91754-2; 978-3-319-91755-9},
 year = {2018},
 publisher = {Cham: Springer},
 language = {English},
 doi = {10.1007/978-3-319-91755-9},
 zbMATH = {6897812},
 Zbl = {1409.53001}
}

@phdthesis{Li26,
	author  = {Zhenhao Li},
	title   = {Optimal transport and (dynamic) synthetic curvature bounds on metric measure spaces},
	school  = {Universität Bielefeld},
	address = {Bielefeld},
	year    = {2026},
	doi     = {10.4119/unibi/3014278},
}

@Article{LottVillani09,
 Author = {Lott, John and Villani, C\'{e}dric},
 Title = {Ricci curvature for metric-measure spaces via optimal transport},
 FJournal = {Annals of Mathematics. Second Series},
 Journal = {Ann. Math. (2)},
 ISSN = {0003-486X},
 Volume = {169},
 Number = {3},
 Pages = {903--991},
 Year = {2009},
 Language = {English},
 DOI = {10.4007/annals.2009.169.903},
 URL = {annals.math.princeton.edu/annals/2009/169-3/p04.xhtml},
 zbMATH = {5578758},
 Zbl = {1178.53038}
}

@article {Mccann2001,
    AUTHOR = {McCann, Robert J.},
     TITLE = {Polar factorization of maps on {R}iemannian manifolds},
   JOURNAL = {Geom. Funct. Anal.},
  FJOURNAL = {Geometric and Functional Analysis},
    VOLUME = {11},
      YEAR = {2001},
    NUMBER = {3},
     PAGES = {589--608},
      ISSN = {1016-443X},
   MRCLASS = {58E15 (46N10 49Q20 53C20)},
  MRNUMBER = {1844080},
MRREVIEWER = {Lucio Renato Berrone},
       DOI = {10.1007/PL00001679},
       URL = {https://doi.org/10.1007/PL00001679},
}

@article {McCann-Topping,
    AUTHOR = {McCann, Robert J. and Topping, Peter M.},
     TITLE = {Ricci flow, entropy and optimal transportation},
   JOURNAL = {Amer. J. Math.},
  FJOURNAL = {American Journal of Mathematics},
    VOLUME = {132},
      YEAR = {2010},
    NUMBER = {3},
     PAGES = {711--730},
      ISSN = {0002-9327},
   MRCLASS = {53C44 (49Q20)},
  MRNUMBER = {2666905},
MRREVIEWER = {Sara Daneri},
       DOI = {10.1353/ajm.0.0110},
       URL = {https://doi.org/10.1353/ajm.0.0110},
}

@article{Schick01,
	author = {Schick, Thomas},
	title = {Manifolds with boundary and of bounded geometry.},
	fjournal = {Mathematische Nachrichten},
	journal = {Math. Nachr.},
	issn = {0025-584X},
	volume = {223},
	pages = {103--120},
	year = {2001},
	language = {English},
	doi = {10.1002/1522-2616(200103)223:1<103::AID-MANA103>3.0.CO;2-S},
	zbMATH = {1606240},
	Zbl = {1037.53016}
}

@article{Shenfeld24,
	author = {Shenfeld, Yair},
	title = {Matrix displacement convexity along density flows},
	fjournal = {Archive for Rational Mechanics and Analysis},
	journal = {Arch. Ration. Mech. Anal.},
	issn = {0003-9527},
	volume = {248},
	number = {5},
	pages = {41},
	note = {Id/No 74},
	year = {2024},
	language = {English},
	doi = {10.1007/s00205-024-02021-8},
	zbMATH = {7922985},
	Zbl = {1554.35340}
}

@incollection{Sparks11,
	author = {Sparks, James},
	title = {Sasaki-{Einstein} manifolds},
	booktitle = {Geometry of special holonomy and related topics},
	isbn = {978-1-57146-211-4},
	pages = {265--324},
	year = {2011},
	publisher = {Somerville, MA: International Press},
	language = {English},
	zbMATH = {6085995},
	Zbl = {1256.53037}
}

@article{Stenzel,
	author = {Stenzel, Matthew B.},
	title = {Ricci-flat metrics on the complexification of a compact rank one symmetric space},
	fjournal = {Manuscripta Mathematica},
	journal = {Manuscr. Math.},
	issn = {0025-2611},
	volume = {80},
	number = {2},
	pages = {151--163},
	year = {1993},
	language = {English},
	doi = {10.1007/BF03026543},
	url = {https://eudml.org/doc/155865},
	zbMATH = {492797},
	Zbl = {0811.53049}
}

@incollection{Sturm03,
	author = {Sturm, Karl-Theodor},
	title = {Probability measures on metric spaces of nonpositive curvature},
	isbn = {0-8218-3383-9},
	year = {2003},
	publisher = {Providence, RI: American Mathematical Society (AMS)},
	language = {English},
	zbMATH = {2051142},
	Zbl = {1040.60002}
}

@article {Sturm05,
    AUTHOR = {Sturm, Karl-Theodor},
     TITLE = {Convex functionals of probability measures and nonlinear
              diffusions on manifolds},
   JOURNAL = {J. Math. Pures Appl. (9)},
  FJOURNAL = {Journal de Math\'{e}matiques Pures et Appliqu\'{e}es. Neuvi\`eme S\'{e}rie},
    VOLUME = {84},
      YEAR = {2005},
    NUMBER = {2},
     PAGES = {149--168},
      ISSN = {0021-7824},
   MRCLASS = {49Q20 (28A33 60B05)},
  MRNUMBER = {2118836},
MRREVIEWER = {Erik J. Balder},
       DOI = {10.1016/j.matpur.2004.11.002},
       URL = {https://doi.org/10.1016/j.matpur.2004.11.002},
}

@Article{sturm2006--1,
 Author = {Sturm, Karl-Theodor},
 Title = {On the geometry of metric measure spaces. {I}},
 FJournal = {Acta Mathematica},
 Journal = {Acta Math.},
 ISSN = {0001-5962},
 Volume = {196},
 Number = {1},
 Pages = {65--131},
 Year = {2006},
 Language = {English},
 DOI = {10.1007/s11511-006-0002-8},
 zbMATH = {5049051},
 Zbl = {1105.53035}
}

@Article{Sturm06-2,
 Author = {Sturm, Karl-Theodor},
 Title = {On the geometry of metric measure spaces. {II}},
 FJournal = {Acta Mathematica},
 Journal = {Acta Math.},
 ISSN = {0001-5962},
 Volume = {196},
 Number = {1},
 Pages = {133--177},
 Year = {2006},
 Language = {English},
 DOI = {10.1007/s11511-006-0003-7},
 zbMATH = {5049052},
 Zbl = {1106.53032}
}

@article{Sturm2017remarks,
 AUTHOR = {Sturm, Karl-Theodor},
     TITLE = {Remarks about synthetic upper {R}icci bounds for metric
              measure spaces},
   JOURNAL = {Tohoku Math. J. (2)},
  FJOURNAL = {The Tohoku Mathematical Journal. Second Series},
    VOLUME = {73},
      YEAR = {2021},
    NUMBER = {4},
     PAGES = {539--564},
      ISSN = {0040-8735},
   MRCLASS = {53C23 (49Q20)},
  MRNUMBER = {4355059},
MRREVIEWER = {Wei Zhao},
       DOI = {10.2748/tmj.20200717},
       URL = {https://doi.org/10.2748/tmj.20200717},
}

@article {Sturm2018Super,
    AUTHOR = {Sturm, Karl-Theodor},
     TITLE = {Super-{R}icci flows for metric measure spaces},
   JOURNAL = {J. Funct. Anal.},
  FJOURNAL = {Journal of Functional Analysis},
    VOLUME = {275},
      YEAR = {2018},
    NUMBER = {12},
     PAGES = {3504--3569},
      ISSN = {0022-1236},
   MRCLASS = {53C44 (31E05 37D40 53C23)},
  MRNUMBER = {3864508},
MRREVIEWER = {Qu\cfac{o}c Anh Ng\^{o}},
       DOI = {10.1016/j.jfa.2018.07.014},
       URL = {https://doi.org/10.1016/j.jfa.2018.07.014},
}

@Article{vonRenesseSturm,
 Author = {von Renesse, Max-K. and Sturm, Karl-Theodor},
 Title = {Transport inequalities, gradient estimates, entropy and {Ricci} curvature},
 FJournal = {Communications on Pure and Applied Mathematics},
 Journal = {Commun. Pure Appl. Math.},
 ISSN = {0010-3640},
 Volume = {58},
 Number = {7},
 Pages = {923--940},
 Year = {2005},
 Language = {English},
 DOI = {10.1002/cpa.20060},
 zbMATH = {2183962},
 Zbl = {1078.53028}
}

@book {Villani,
    AUTHOR = {Villani, C\'{e}dric},
     TITLE = {Optimal transport},
    SERIES = {Grundlehren der Mathematischen Wissenschaften},
    VOLUME = {338},
      NOTE = {Old and new},
 PUBLISHER = {Springer-Verlag, Berlin},
      YEAR = {2009},
     PAGES = {xxii+973},
      ISBN = {978-3-540-71049-3},
   MRCLASS = {49-02 (28A75 37J50 49Q20 53C23 58E30)},
  MRNUMBER = {2459454},
MRREVIEWER = {Dario Cordero-Erausquin},
       DOI = {10.1007/978-3-540-71050-9},
       URL = {https://doi.org/10.1007/978-3-540-71050-9},
}

@book{Wang_book14,
	author = {Wang, Feng-Yu},
	title = {Analysis for diffusion processes on {Riemannian} manifolds},
	fseries = {Advanced Series on Statistical Science \& Applied Probability},
	series = {Adv. Ser. Stat. Sci. Appl. Probab.},
	issn = {1793-091X},
	volume = {18},
	isbn = {978-981-4452-64-9; 978-981-4452-66-3},
	year = {2014},
	publisher = {Hackensack, NJ: World Scientific},
	language = {English},
	doi = {10.1142/8737},
	zbMATH = {6198040},
	Zbl = {1296.58001}
}
\end{document}